\documentclass{amsart}

\pdfmapfile{=tengwarscript.map}
\usepackage[annatar]{tengwarscript}

\usepackage{ragged2e}
\usepackage{a4wide}
\usepackage[utf8]{inputenc}
\usepackage{enumerate}

\usepackage{amsmath}
\usepackage{amsfonts}
\usepackage{mathabx}

\usepackage{bbm}
\usepackage{mathrsfs}
\usepackage{amsthm}
\usepackage{verbatim}
\usepackage{upgreek}
\usepackage{relsize}
\usepackage{dsfont}
\usepackage{graphicx}
\usepackage{amssymb}
\usepackage{xcolor}

\newcommand{\eps}{\varepsilon}
\newcommand{\ov}{\overline}
\newcommand{\id}{\textnormal{id}}
\newcommand{\mc}{\mathcal}

\newcommand{\mrm}{\mathrm}

\newcommand{\mf}{\mathfrak}
\newcommand{\msf}{\mathsf}
\newcommand{\I}{\mathbbm{1}}

\newcommand{\vp}{\varphi}

\newcommand{\md}{\operatorname{d}\!}
\newcommand{\op}{\operatorname{op}}
\newcommand{\cst}{\ifmmode \mathrm{C}^* \else $\mathrm{C}^*$\fi}

\newcommand{\la}{\langle}
\newcommand{\ra}{\rangle}

\newcommand{\NN}{\mathbb{N}}
\newcommand{\RR}{\mathbb{R}}

\newcommand{\CC}{\mathbb{C}}
\newcommand{\ZZ}{\mathbb{Z}}
\newcommand{\GG}{\mathbb{G}}
\newcommand{\HH}{\mathbb{H}}

\newcommand{\vv}{\mathrm{V}}
\newcommand{\Vv}{\mathds{V}}
\newcommand{\vV}{\text{\reflectbox{$\Vv$}}\:\!}

\newcommand{\ww}{\mathrm{W}}
\newcommand{\WW}{{\mathds{V}\!\!\text{\reflectbox{$\mathds{V}$}}}}
\newcommand{\Ww}{\mathds{W}}
\newcommand{\wW}{\text{\reflectbox{$\Ww$}}\:\!}

\newcommand{\wot}{\ifmmode \textsc{wot} \else \textsc{wot}\fi}
\newcommand{\sot}{\ifmmode \textsc{sot} \else \textsc{sot}\fi}
\newcommand{\sots}{\ifmmode \textsc{sot}^* \else \textsc{sot}$^*$\fi}
\newcommand{\ssot}{\ifmmode \sigma\textsc{-sot} \else $\sigma$-\textsc{sot }\fi}
\newcommand{\ssots}{\ifmmode \sigma\textsc{-sot}^* \else $\sigma$-\textsc{sot }$^*$\fi}
\newcommand{\swot}{\ifmmode \sigma\textsc{-wot} \else $\sigma$-\textsc{wot}\fi}

\newcommand{\wh}{\widehat}
\newcommand{\wc}{\widecheck}
\newcommand{\wt}{\widetilde}
\newcommand{\whG}{\widehat{\GG}}

\newcommand{\oon}{\operatorname}

\renewcommand{\restriction}{\mathord{\upharpoonright}}
\newcommand{\rest}{\restriction}

\DeclareMathOperator{\lin}{span}

\DeclareMathOperator{\Tr}{Tr}
\DeclareMathOperator{\B}{B}

\DeclareMathOperator{\M}{M}

\DeclareMathOperator{\N}{N}
\DeclareMathOperator{\OnP}{P}
\DeclareMathOperator{\Dom}{Dom}

\DeclareMathOperator{\LL}{L}
\DeclareMathOperator{\HS}{HS}

\DeclareMathOperator{\CB}{CB}

\DeclareMathOperator{\Aut}{Aut}
\DeclareMathOperator{\Ad}{Ad}
\DeclareMathOperator{\ext}{ext}
\DeclareMathOperator{\nsf}{nsf}
\DeclareMathOperator{\vN}{vN}

\newtheorem{theorem}{Theorem}[section]

\newtheorem{proposition}[theorem]{Proposition}
\newtheorem{lemma}[theorem]{Lemma}
\theoremstyle{definition}
\newtheorem{corollary}[theorem]{Corollary}
\newtheorem{remark}[theorem]{Remark}

\newtheorem{definition}[theorem]{Definition}
\newtheorem{example}[theorem]{Example}

\numberwithin{equation}{section}

\begin{document}
\title[Actions of quantum groups]{
Actions of quantum groups: the standard implementation and relatively invariant weights}

%

\author{Jacek Krajczok}
\address{Institute of Mathematics of the Polish Academy of Sciences, Śniadeckich 8, 00-656 Warsaw, Poland}
\email{jkrajczok@impan.pl}

\thanks{}





\date{}

\begin{abstract}
We establish a characterisation of the standard implementation for actions of locally compact quantum groups on von Neumann algebras, extending results of Daws-Skalski-Viselter and Haagerup. As a consequence, we prove that taking the standard implementation respects morphisms of quantum groups. We also prove that several natural conditions expressing relative $\rho$-invariance for n.s.f.~weights are equivalent.
\end{abstract}

\maketitle

\section{Introduction}
Any von Neumann algebra can be represented faithfully on a Hilbert space in a number of different ways. While some of the properties do not depend on the choice of representation, others -- such as e.g.~weak operator topology -- do. Luckily, any von Neumann algebra admits a canonical and unique representation, in a certain precise sense. Haagerup in \cite{HaagerupStandard} introduced the notion of the \emph{standard representation} (or \emph{standard form}) of a von Neumann algebra $\N$. It is a $4$-tuple $(\msf{H},\pi,J,\mf{P})$ consisting of a Hilbert space $\msf{H}$, faithful, unital, normal representation $\pi\colon \N\rightarrow \B(\msf{H})$, self-adjoint anti-unitary $J$ on $\msf{H}$ and cone $\mf{P}\subseteq \msf{H}$, satisfying certain conditions. Crucially, this structure exists for any von Neumann algebra, and two such representations are conjugated by a unique unitary which respects $J$ and $\mf{P}$. Consequently, standard representation can be thought of as a structure that is inherent to $\N$ itself. Furthermore, it has several convenient properties: any normal functional on $\N$ is a vector functional in this representation, and  $J \pi(\N) J= \pi(\N)'$ -- compare these statements to e.g.~the defining representation of $\B(\msf{K})$ on $\msf{K}$.

Any automorphism $\alpha\in \Aut(\N)$ gives rise to a new representation $\pi \alpha\colon \N\rightarrow \B(\msf{H})$. The uniqueness of standard representation implies that there is a unique unitary $u\in \B(\msf{H})$ which satisfies $\pi(\alpha(n))=u \pi(n) u^*\,(n\in \N)$ and $uJ=Ju, u\mf{P}=\mf{P}$. In other words, in the standard representation any automorphism has a unique unitary implementation which respects \emph{modular conjugation} $J$ and \emph{positive cone} $\mf{P}$. Note that we do not impose the existence of a n.s.f.~weight on $\N$ which is invariant under $\alpha$; such a weight might fail to exist.

More generally, if $\alpha\colon G\curvearrowright \N$ is a weak$^*$-continuous action of a locally compact group $G$ on $\N$, then each automorphism $\alpha_g\,(g\in G)$ is implemented by the corresponding unitary $u_g\in \B(\msf{H})$, and the resulting family $(u_g)_{g\in G}$ forms a \sot-continuous representation of $G$; the \emph{standard implementation} of action $\alpha$ (\cite[Corollary 3.6]{HaagerupStandard}, see also Example \ref{ex5}). Haagerup's results show existence of the representation $(u_g)_{g\in G}$ implementing $\alpha$ and provide conditions $u_g J =J u_g, u_g\mf{P}=\mf{P}$ which characterise this implementation uniquely. The first part of this manuscript is concerned with extending such uniqueness statement to the quantum setting, which we explain next.
\\

A far reaching generalisation of the theory of locally compact groups was proposed by Kustermans and Vaes \cite{KustermansVaes1, KustermansVaes}, see also \cite{VanDaeleVN}. Objects in this theory are called \emph{locally compact quantum groups}. Every such $\GG$ is described by a von Neumann algebra $\LL^{\infty}(\GG)$, comultiplication $\Delta\colon \LL^{\infty}(\GG)\rightarrow \LL^{\infty}(\GG)\bar\otimes \LL^{\infty}(\GG)$ and left (resp.~right) invariant n.s.f.~weight $\vp,\psi$ on $\LL^{\infty}(\GG)$; the Haar integrals. The class of locally compact quantum groups contains classical theory of locally compact groups and is closed under (the generalisation of) Pontryagin duality. Interesting examples of quantum groups have been constructed, both in compact and non-compact case, see \cite{Timmermann, NeshveyevTuset, VoigtYuncken}.

Similarly as their classical counterparts, locally compact quantum groups can act on von Neumann algebras (see Preliminaries for further explanation). For an action $\alpha\colon \GG\curvearrowright \N$ and a fixed arbitrary n.s.f.~weight $\theta$ on $\N$,
 Vaes in \cite{VaesUnitary} constructed a \emph{standard implementation} of $\alpha$: a particular unitary representation $\mrm{U}^\alpha\in \LL^{\infty}(\GG)\bar\otimes \B(\msf{H}_\theta)$ of $\GG$, which implements $\alpha$ and respects modular conjugations:
\begin{enumerate}
\item[(1)] $\alpha(n)=\mrm{U}^{\alpha *}(\I\otimes n)\mrm{U}^{\alpha}\quad(n\in\N)$,
\item[(2)] $(J_{\wh\vp}\otimes J_\theta)\mrm{U}^{\alpha *} (J_{\wh\vp}\otimes J_\theta)=\mrm{U}^{\alpha}$.
\end{enumerate}
Furthermore, Daws-Skalski-Viselter showed in \cite{DawsSkalskiViselter} that
\begin{enumerate}
\item[(3)] $\mrm{U}^{\alpha}$ respects the positive cone
\end{enumerate}
(see \cite[Proposition 4.10]{DawsSkalskiViselter} and also Remark \ref{rem3} for a precise condition and discussion). A natural question arises: do conditions $(1),\, (2),\, (3)$ characterise $\mrm{U}^\alpha$? The authors of \cite{DawsSkalskiViselter} established this claim for discrete quantum groups, and we prove that it is the case in full generality, using different reasoning (Theorem \ref{thm1}).

As an important corollary of this characterisation we show that the standard implementation of the action given by a composition with morphism, $\GG\rightarrow \HH\curvearrowright \N$, coincides with applying the morphism $\GG\rightarrow \HH$ to the standard implementation of $\HH\curvearrowright \N$, namely $\mathds{U}^{\alpha_{\GG}}=( \Pi\otimes\id)(\mathds{U}^{\alpha _{\HH}})$ (see Corollary \ref{cor1} and the preceding paragraph). This result is motivated by a work of Kenny De Commer and the author on braided locally compact quantum groups \cite{DeCommerKrajczokBLCQG}. However, to prove Corollary \ref{cor1} we will first need to establish a universal version of conditions $(1),\, (2),\, (3)$ (Proposition \ref{prop3} and Theorem \ref{thm5}) which will require a considerable effort.\\

The second topic considered in this article is the notion of relative invariance for weights. Assume that $\alpha\colon \GG\curvearrowright \N$ is an action, $\theta$ is a n.s.f.~weight on $\N$ and $\rho$ is a strictly positive, self-adjoint group-like operator affiliated with $\LL^{\infty}(\GG)$. One way of saying that $\theta$ is \emph{relatively $\rho$-invariant} is to impose the condition
\begin{equation}\tag{RI1}
\theta\bigl((\omega_\xi\otimes\id)\alpha(n)\bigr)=\|\rho^{1/2} \xi\|^2 \theta(n)
\end{equation}
for $\xi \in \Dom(\rho^{1/2})$ and integrable $n\in \N^+, \theta(n)<+\infty$. An a priori stronger condition states
\begin{equation}\tag{RI3}
(\id\otimes\id\otimes \theta)
(\id\otimes\alpha)(X)
=
(\id\otimes \theta)(X)\otimes \rho
\end{equation}
for an arbitrary von Neumann algebra $\OnP$ and $X\in (\OnP\bar\otimes\N)^+$. Thus in (RI3) one does not assume integrability of the element $X$, and the amplification with arbitrary coefficient algebra $\OnP$ is allowed. The first condition (RI1) was introduced in \cite{VaesUnitary}. Such weights often arise in practice, for example the dual weight on crossed product satisfies (RI1). It is also the invariance condition used in the definition of Haar integrals (with $\rho=\I$). We will also introduce conditions (RI2), (RI4) and prove in Theorem \ref{thm2} that all these variations on relative invariance are equivalent.

While the difference between (RI1) and (RI3) looks rather technical, a stronger condition is an important technical tool which is sometimes indispensable. For example, (RI3) for the left Haar integral with respect to the action $\Delta$ (proved by Kustermans and Vaes in \cite[Proposition 3.1]{KustermansVaes}) was crucially used in the construction of the dual weight on crossed product \cite[Definition 3.1]{VaesUnitary} and the authors of \cite{BrannanChirvasituViselter} used (RI2) in their work on relative invariant weights on homogeneous spaces $\LL^{\infty}(\GG/\HH)$ (see Remark \ref{rem2}).\\

We have described a number of examples of the standard implementation in \ref{ex1} -- \ref{ex6}, and examples of relatively invariant weights in \ref{ex8} -- \ref{ex9}.\\

\textbf{Acknowledgements.} 
I am grateful to Kenny De Commer for discussions concerning matters touched in this manuscript. This work was partially supported by FWO grant 1246624N.\\

\textbf{AI statement.} Arguments presented in this work were partially obtained as a result of discussions with ChatGPT 5.6 Sol and ChatGPT 6 Astra, later clarified, cleaned and improved by the author. The author takes full responsibility for the validity of presented proofs.

\section{Preliminaries and notation}

Whenever $\N$ is a von Neumann algebra, we denote by $\oon{W}_{\nsf}(\N)$ the set of n.s.f.~(normal, semifinite, faithful) weights on $\N$. For $\theta\in \oon{W}_{\nsf}(\N)$, $\msf{H}_\theta$ is the GNS Hilbert space with faithful representation $\pi_\theta\colon \N\rightarrow \B(\msf{H}_\theta)$, GNS map $\Lambda_\theta\colon \mf{N}_\theta\rightarrow \msf{H}_\theta$, $J_\theta$ is the modular conjugation and $\mf{P}_\theta$ is the positive cone. The $4$-tuple $(\msf{H}_\theta,\pi_\theta,J_\theta,\mf{P}_\theta)$ is a standard representation of $\N$, but we will typically neglect writing representation $\pi_\theta$. We follow the usual notation $\mf{M}_\theta^+=\{n\in\N^+\mid \theta(n)<+\infty\}$, $\mf{M}_\theta=\lin \mf{M}_\theta^+$ and $\mf{N}_\theta=\{n\in\N\mid \theta(n^*n)<+\infty\}$ \cite{StratilaZsido, Stratila, HaagerupStandard}, and use inner products linear in the second variable.

If $\N\subseteq \B(\msf{H})$ is a von Neumann algebra and $A$ is a closed, densely-defined operator on $\msf{H}$, we denote its domain by $\Dom(A)$ and the affiliation relation by $A\eta \N$ \cite{StratilaZsido}. If $A$ is positive and $\ker(A)=\{0\}$, then we say that $A$ is strictly positive. We will use the notion of \emph{extended positive part} $\N^{+,\ext}$ \cite{HaagerupOVWI, HaagerupOVWII, Stratila}. It consists of all maps $X\colon \N_*^+\rightarrow \left[0,\infty\right]$ which are additive, positively homogeneous and lower semicontinuous; one also writes $\omega(X)=X(\omega)$ for $\omega\in \N_*^+$. Any $X\in \N^{+,\ext}$ has a spectral resolution: there are unique projection $p\in \N$ and positive, self-adjoint $A \eta\, p\N p\subseteq \B(p\msf{H})$ such that
\begin{equation}\label{eq42}
X=A+\infty(\I-p)= \int_0^\infty \lambda \md e_{A}(\lambda)+\infty (\I-p),
\end{equation}
where $e_A$ is the spectral measure of $A$ (see \cite{HaagerupOVWI, SoltanPrimer}). Conversely, any such pair gives rise to an element of $\N^{+,\ext}$, in particular $\N^+\subseteq \N^{+,\ext}$. We can think of elements in $\N^{+,\ext}$ as unbounded positive operators with possibly non-trivial infinite and zero parts.

For $X\in \N^{+,\ext}$ with associated $A,p$, we write $\Dom(X)=\Dom(A)$, and define a closed quadratic form $q_X$ with domain $\Dom(q_X)=\Dom(A^{1/2})$, via $q_X(\xi)=\|A^{1/2}\xi\|^2\,(\xi\in \Dom(q_X))$ \cite{Stratila, Kato}. $q_X$ induces in the standard way a sesquilinear form $q_X\colon \Dom(q_X)^{\times 2}\rightarrow \CC$. We also extend the domain of $q_X$ and write $q_X(\xi)=\infty$ if $\xi\in\msf{H}\setminus \Dom(q_X)$, so that $q_X(\xi)=\omega_\xi(X)$.

One easily checks that normal, unital $*$-homomorphisms can be applied to elements of $\N^{+,\ext}$. We can also define the tensor product $X\otimes Y\in (\N\bar\otimes \M)^{+,\ext}$, where $X=A_X+\infty(\I-p_X)\in \N^{+,\ext}$ and $Y\eta \M$ is strictly positive, self-adjoint. Element $X\otimes Y$ corresponds to the pair $\bigl( p_X\otimes \I,  A_X\otimes Y\bigr)$. One obtains the same result by setting $X\otimes Y=\sup_{n\in\NN} X_n\otimes Y_n$ for $X_n=\int_0^n \lambda \md e_{A_X}(\lambda)+n (\I-p_X), Y_n=\int_0^n\lambda\md e_{Y}(\lambda)$.\\
 
We will work with locally compact quantum groups, as introduced by Kustermans and Vaes \cite{KustermansVaes1,KustermansVaes}. By definition, a locally compact quantum group $\GG$ is a $4$-tuple $(\LL^{\infty}(\GG),\Delta_{\GG},\vp_{\GG},\psi_{\GG})$, where $\LL^{\infty}(\GG)$ is a von Neumann algebra, $\Delta_{\GG}\colon \LL^{\infty}(\GG)\rightarrow \LL^{\infty}(\GG)\bar\otimes \LL^{\infty}(\GG)$ is a normal, injective, unital $*$-homomorphism (\emph{comultiplication}) satisfying $(\Delta_{\GG}\otimes\id)\Delta_{\GG}=(\id\otimes\Delta_{\GG})\Delta_{\GG}$ and $\vp_{\GG},\psi_{\GG}\in \oon{W}_{\nsf}(\LL^{\infty}(\GG))$ are weights satisfying left/right invariance condition ($\vp_{\GG}$ satisfies (RI1) with action $\Delta_{\GG}$ and $\psi_{\GG}$ its right-handed analog). The GNS Hilbert spaces of $\vp_{\GG}$ and $\psi_{\GG}$ can be identified, and will be denoted by $\LL^2(\GG)$. Left and right Haar integral are related by the \emph{modular element} $\delta_{\GG}\eta\LL^{\infty}(\GG)$; it holds that $\psi_{\GG}=(\vp_{\GG})_{\delta_{\GG}}$ \cite{KustermansVaes}. If $\N$ is a von Neumann algebra, an \emph{action} $\alpha\colon \GG\curvearrowright \N$ is a unital, normal, injective $*$-homomorphism $\alpha\colon \N\rightarrow \LL^{\infty}(\GG)\bar\otimes \N$ satisfying $(\Delta_{\GG}\otimes\id)\alpha=(\id\otimes\alpha)\alpha$. A \emph{unitary representation} of $\GG$ on Hilbert space $\msf{H}$ is a unitary $\mrm{U}\in \LL^{\infty}(\GG)\bar\otimes \B(\msf{H})$ satisfying $(\Delta_{\GG}\otimes\id)\mrm{U}=\mrm{U}_{[13]} \mrm{U}_{[23]}$; we will consider only unitary representations. If $\N\subseteq \B(\msf{H})$ and $\alpha(n)=\mrm{U}^*(\I\otimes n)\mrm{U}$, then we say that representation $\mrm{U}$ implements action $\alpha$. Note that we primarily use \emph{left} actions and \emph{left} representations, but all our results have obvious \emph{right} counterparts. Of particular importance is the \emph{left regular representation} (known also as \emph{Kac-Takesaki operator}) $\ww^{\GG}$, which is a representation of $\GG$ on $\LL^2(\GG)$. It implements comultiplication in the following way: $\Delta_{\GG}(x)=\ww^{\GG *}(\I\otimes x)\ww^{\GG}\,(x\in\LL^{\infty}(\GG))$. Haar integrals $\vp_{\GG},\psi_{\GG}$ give rise to \emph{modular automorphism group} $\sigma^{\vp_{\GG}}_t,\sigma^{\psi_{\GG}}_t\in \Aut(\LL^{\infty}(\GG))\,(t\in\RR)$. There is also a \emph{unitary antipode} $R_{\GG}$ on $\LL^{\infty}(\GG)$ and a \emph{scaling group} $\tau^{\GG}_t\in\Aut(\LL^{\infty}(\GG))\,(t\in\RR)$, implemented by $R_{\GG}(x)=J_{\vp_{\whG}} x^* J_{\vp_{\whG}}, \tau^{\GG}_t(x)=P_{\GG}^{ it} x P_{\GG}^{ -it}\,(x\in \LL^{\infty}(\GG),t\in\RR)$, where $P_{\GG}$ is a strictly positive, self-adjoint operator on $\LL^2(\GG)$, self-dual in the sense $P_{\whG}=P_{\GG}$. In fact, $P_{\GG}^{it}$ is the standard implementation of $\tau^{\GG}_t$, see Example \ref{ex6}. Space $\LL^{\infty}(\whG)=\{(\omega\otimes\id)\ww^{\GG}\mid \omega\in \LL^1(\GG)\}^{-\swot}$ is a von Neumann algebra, and as the notation suggests, it is a von Neumann algebra of (measurable, bounded) functions on quantum group $\whG$, which is the Pontryagin dual of $\GG$. Its comultiplication is implemented by $\ww^{\whG}=\ww^{\GG *}_{[21]}$. There is also a \emph{right regular representation} $\vv^{\GG}$ which implements comultiplication via $\Delta_{\GG}(x)=\vv^{\GG}(x\otimes\I)\vv^{\GG *}$. It is a \emph{right unitary representation} in the sense that $\vv^{\GG}\in \B(\LL^2(\GG))\bar\otimes \LL^{\infty}(\GG)$ is unitary and $(\id\otimes \Delta_{\GG})\vv^{\GG}=\vv^{\GG}_{[12]}\vv^{\GG}_{[13]}$. We also have $\ww^{\GG}\in \LL^{\infty}(\GG)\bar\otimes \LL^{\infty}(\whG)$, $\vv^{\GG}\in \LL^{\infty}(\wc\GG)\bar\otimes \LL^{\infty}(\GG)$ and $(\id\otimes\Delta_{\whG})(\ww^{\GG})=\ww^{\GG}_{[13]} \ww^{\GG}_{[12]}$, $(\Delta_{\wc\GG}\otimes\id)(\vv^{\GG})=\vv^{\GG}_{[13]} \vv^{\GG}_{[23]}$, where $\wc\GG=\whG'$ \cite[Section 4]{KustermansVaes}. These objects satisfy a plethora of relations, for which we refer to \cite{KustermansVaes, VanDaeleVN, KrajczokTypeI, Kustermans}.

One can also study quantum groups in \cst-algebraic language; both approaches are equivalent and often one needs to use both. The space $\mrm{C}_0(\GG)=\{(\id\otimes\omega)\ww^{\GG}\mid \omega\in \LL^1(\whG)\}^{-}$ is a weak$^*$-dense \cst-subalgebra of $\LL^{\infty}(\GG)$, which plays the role of continuous functions on $\GG$. It admits a universal counterpart $\mrm{C}_0^u(\GG)$ \cite{Kustermans}, which comes equipped with canonical surjective \emph{reducing morphism} $\pi_{\GG}\colon \mrm{C}_0^u(\GG)\rightarrow \mrm{C}_0(\GG)$. A lot of structure (such as e.g.~comultiplication, modular automorphism group, scaling group and unitary antipode) restricts to $\mrm{C}_0(\GG)$ and lifts to $\mrm{C}_0^u(\GG)$. For example, $\Delta_{\GG}$ restricts to a non-degenerate $*$-homomorphism $\mrm{C}_0(\GG)\rightarrow \M(\mrm{C}_0(\GG)\otimes \mrm{C}_0(\GG))$ and Kac-Takesaki operator $\ww^{\GG}$ admits lift to a unitary $\wW^{\GG}\in \M(\mrm{C}_0(\GG)\otimes \mrm{C}_0^u(\whG))$ satisfying $(\id\otimes \pi_{\whG})(\wW^{\GG})=\ww^{\GG}$. There are also versions with left or both legs lifted $\Ww^{\GG}, \WW^{\GG}$; these unitaries have analogous properties with respect to comultiplication as $\ww^{\GG}$. Using them, one can define also half-lifted comultiplications $\Delta_{\GG}^{u,r}\colon \mrm{C}_0(\GG)\rightarrow \M(\mrm{C}_0^u(\GG)\otimes \mrm{C}_0(\GG))$ and similarly $\Delta_{\GG}^{r,u}$.

\cst-algebra $\mrm{C}_0^u(\whG)$ plays the role of full group \cst-algebra of $\GG$, and $\wW^{\GG}$ the role of universal representation: whenever $\mathrm{U}$ is a representation of $\GG$ on $\msf{H}$, there is a unique non-degenerate representation $\phi\colon \mrm{C}_0^u(\whG)\rightarrow \B(\msf{H})$ such that $\mathrm{U}=(\id\otimes \phi)(\wW^{\GG})$. Using this, we can consider lift $\mathds{U}=(\id\otimes\phi)(\WW^{\GG})$. It satisfies $(\pi_{\GG}\otimes\id)(\mathds{U})=\mathrm{U}$, $(\Delta_{\GG}^u\otimes \id)(\mathds{U})=\mathds{U}_{[13]} \mathds{U}_{[23]}$  and is the unique unitary in $\M(\mrm{C}_0^u(\GG)\otimes \mathcal{K}(\msf{H}))$ with these properties \cite{Kustermans}. Action $\alpha\colon \GG\curvearrowright \N$ admits a universal lift: $\alpha^{u}\colon \N\rightarrow \mrm{C}_0^u(\GG)^{**}\bar\otimes \N$, which has similar properties as $\alpha$ with respect to (half) lifted von Neumann algebraic comultiplication \cite[Section 4.3]{DeCommerKrajczok}. If it does not lead to confusion, we will ease the notation and neglect writing (sub/super)-scripts, e.g.~$\Delta=\Delta_{\GG}$. \\

Whenever $A$ is a \cst-algebra, its bidual $A^{**}$ carries a structure of a von Neumann algebra and $A\subseteq A^{**}$ is weak$^*$-dense \cite{StratilaZsido}. We have the following useful property: if $\M$ is a von Neumann algebra, and $\pi\colon A\rightarrow \M$ is a non-degenerate $*$-homomorphism, then there is a unique extension to a normal, unital $*$-homomorphism $\pi^{\vN}\colon A^{**}\rightarrow \M$. 

We will also need, in a limited scope, the notion of an element affiliated with a \cst-algebra \cite{WoronowiczUnbounded, Kustermans1997functional}, still denoted $a\eta A$. In particular, we will use their functional calculus and the fact that one can apply non-degenerate $*$-homomorphism to such elements. Let us also recall the following, well-known construction. Let $\rho$ be a strictly positive, self-adjoint operator affiliated with $\LL^{\infty}(\GG)$, which is group-like in the sense $\Delta_{\GG}(\rho)=\rho\otimes \rho$. Then for $t\in \RR$, $\rho^{it}$ is a group-like unitary in $\M(\mrm{C}_0(\GG))$, hence it lifts uniquely to a group-like unitary in $\M(\mrm{C}_0^u(\GG))$ (it can be seen as a one dimensional representation), which by foreshadowing we denote $\rho^{u it}\in \M(\mrm{C}_0^u(\GG))$. One easily checks that the group $(\rho^{u it})_{t\in\RR}$ is strictly-continuous, hence it indeed admits a strictly positive, self-adjoint generator $\rho^u\,\eta\,\mrm{C}_0^u(\GG)$: combine \cst-algebraic version of Stone's theorem \cite{StoneTheoremCS} with functional calculus \cite[Corollary 7.15, Proposition 7.11]{Kustermans1997functional}. In particular, if $\pi\colon \mrm{C}_0^u(\GG)\rightarrow \B(\msf{K})$ is a non-degenerate representation, then $\pi(\rho^{u})$ is the unique strictly positive, self-adjoint operator on $\msf{K}$ satisfying $\pi(\rho^u)^{it}=\pi(\rho^{u it})\,(t\in\RR)$.

\section{Characterisation of the standard implementation of action}\label{sec:characterisation}
In this section we establish a characterisation of the standard implementation for an action of a LCQG \cite{VaesUnitary}, extending the result of \cite{DawsSkalskiViselter} (see Theorem \ref{thm1} and Remark \ref{rem3}). Afterwards we discuss a number of natural examples of actions and their standard implementations.

Let $\GG$ be a locally compact quantum group with a (left) action on a von Neumann algebra $\alpha\colon \GG\curvearrowright \N$. Fix a n.s.f.~weight $\theta\in \oon{W}_{\nsf}(\N)$; this gives a standard structure $(\msf{H}_\theta,J_\theta,\mf{P}_\theta)$ for $\N$. Then one defines the crossed product
\[
\GG\ltimes \N=\ov{\lin}^{\,\swot}\bigl\{(\hat{x}\otimes\I)\alpha(n)\mid\hat{x}\in \LL^{\infty}(\whG),n\in\N\bigr\}\subseteq \B(\LL^2(\GG)\otimes\msf{H}_\theta),
\]
which admits the dual weight $\wh{\theta}\in \oon{W}_{\nsf}(\GG\ltimes \N)$ \cite[Definition 3.1]{VaesUnitary}. The GNS representation of $\wh{\theta}$ is realised in $\LL^2(\GG)\otimes \msf{H}_\theta$ and $\GG\ltimes \N$ is represented in a standard way on this space: $(\LL^2(\GG)\otimes\msf{H}_\theta,J_{\wh{\theta}},\mf{P}_{\wh\theta})$. Furthermore,
\[
\mrm{U}^\alpha=(J_{\wh\vp}\otimes J_{\theta})J_{\wh{\theta}}
\]
is\footnote{Unitary $\mrm{U}^\alpha$ is the adjoint of the one from \cite{VaesUnitary}; we introduce this change so that $\mrm{U}^\alpha$ is a representation of $\GG$.} a unitary representation of $\GG$ implementing $\alpha$ in the sense that $\alpha(n)=\mrm{U}^{\alpha *}(\I\otimes n)\mrm{U}^\alpha\,(n\in\N)$. Representation $\mrm{U}^{\alpha}$ is called the \emph{standard} implementation \cite[Definition 3.6]{VaesUnitary}.

\begin{remark}\label{rem1}
Strictly speaking, $\mrm{U}^\alpha$ depends on the choice of $\theta$ -- however, change of the weight changes the resulting representation in a convenient way. Indeed, given $\theta_1,\theta_2\in \oon{W}_{\nsf}(\N)$, there is a unique unitary $u$ which conjugates two standard representations of $\N$ \cite[Theorem 2.3]{HaagerupStandard}. This unitary transports between two standard implementations $\mrm{U}^{\alpha,1},\mrm{U}^{\alpha,2}$ (constructed with respect to $\theta_1$ and $\theta_2$) in the sense that $(\I\otimes u)\mrm{U}^{\alpha,1}=\mrm{U}^{\alpha,2}(\I\otimes u)$ \cite[Proposition 4.1]{VaesUnitary}. This reasoning shows also that $\mrm{U}^\alpha$ depends only on the standard representation $(\msf{H}_\theta,J_\theta,\mf{P}_\theta)$ rather than the weight $\theta$ itself.
\end{remark}

In what follows, we will need a particular standard representation of $\M_n\bar\otimes \N$. The natural choice is again a GNS representation: $\Tr_n\otimes\, \theta$ is a n.s.f.~weight on $\M_n\bar\otimes \N$ and its GNS representation is realised in $\HS(\CC^n)\otimes \msf{H}_\theta$ (elements of which will be identified with $\msf{H}_\theta$-valued matrices whenever convenient).

\begin{lemma}\label{lemma3}
We have
\[\begin{split}
\mf{P}_{\Tr_n\otimes\theta}&=
\bigl\{
\sum_{\lambda=1}^{\Lambda} \sum_{i,j=1}^{n} e_{i,j}\otimes x^\lambda_i J_\theta \Lambda_\theta(x^\lambda_j) \mid \Lambda\in \NN, x^\lambda_{1},\dotsc,x^\lambda_n \in \mf{N}_\theta \}^{-}\\
&=
\bigl\{
\sum_{\lambda=1}^{\Lambda} \sum_{i,j=1}^{n} e_{i,j}\otimes x^\lambda_i J_\theta x^\lambda_j J_\theta \zeta^\lambda  \mid \Lambda\in \NN, x^\lambda_{1},\dotsc,x^\lambda_n \in \mf{N}_\theta,\zeta^\lambda\in \mf{P}_\theta \}^{-},
\end{split}\]
where $\{e_{i,j}\}_{i,j=1}^{n}$ are the matrix units in $\M_n$.
\end{lemma}

\begin{proof}
The first expression follows easily from the formula $\mf{P}_{\Tr_n\otimes\,\theta}=\{x J_{\Tr_n\otimes\,\theta}\Lambda_{\Tr_n\otimes\theta}(x)\mid x\in \mf{N}_{\Tr_n\otimes\theta}\}^{-}$ \cite[Page 29]{Stratila}. For the second, observe that
\[
x_i J_\theta x_jJ_\theta \zeta=
x_i J_\theta x_jJ_\theta yJ_\theta \Lambda_\theta(y)=
x_i y J_\theta \Lambda_\theta(x_jy)
\]
for $x_1,\dotsc,x_n\in \mf{N}_\theta$ and $\zeta=yJ_\theta \Lambda_\theta(y)\in \mf{P}_\theta$ and
\[
\sum_{i,j=1}^{n} e_{i,j}\otimes x_iJ_\theta\Lambda_\theta(x_j)=\!
\lim_{\mu\in \Sigma}\sum_{i,j=1}^{n}\! e_{i,j}\otimes x_i f_\mu \sigma^\theta_{i/2}(f_\mu)^* J_\theta\Lambda_\theta(x_j)=\!
\lim_{\mu\in \Sigma}\sum_{i,j=1}^{n} \!e_{i,j}\otimes x_i J_\theta x_j J_\theta \bigl(f_\mu J_\theta \Lambda_\theta(f_\mu)\bigr)
\] where $(f_\mu)_{\mu\in \Sigma}$ is an approximate unit from \cite[Proposition 2.16]{Stratila}.
\end{proof}

The following definition appeared in \cite[Definition 4.8]{DawsSkalskiViselter}, following \cite[Section 2.5]{Sauvageot} in the classical case.

\begin{definition}\label{def1}
Let $\N$ be a von Neumann algebra with a n.s.f.~weight $\theta$, and $\mrm{U}=(\mrm{U}_{i,j})_{i,j=1}^{n}\in \M_n(\B(\msf{H}_\theta))$ a matrix of operators. We say that $\mrm{U}$ \emph{preserves the positive cone} if for every $(\zeta_{i,j})_{i,j=1}^{n} \in \mf{P}_{ \Tr_n\otimes \theta}$, the matrix $(\mrm{U}_{i,j}\zeta_{i,j})_{i,j=1}^{n}$ belongs to the cone $\mf{P}_{\Tr_n\otimes \theta}$.
\end{definition}

Let us stress that in Definition \ref{def1}, the Schur (rather than matrix) multiplication is used. The following is the main result of this section.

\begin{theorem}\label{thm1}
Let $\N$ be a von Neumann algebra with a weight $\theta\in \oon{W}_{\nsf}(\N)$ and action $\alpha\colon \GG\curvearrowright \N$. Let $\mrm{U}\in \LL^{\infty}(\GG)\bar\otimes \B(\msf{H}_\theta)$ be a unitary representation of $\GG$ and $\mrm{U}^\alpha$ the standard implementation of $\alpha$. Then $\mrm{U}=\mrm{U}^\alpha$ if, and only if
\begin{enumerate}
\item $\mrm{U}^*(\I\otimes y)\mrm{U}=\alpha(y)$ for $y\in\N$,
\item $(J_{\wh\vp}\otimes J_\theta)\mrm{U}^*(J_{\wh\vp}\otimes J_\theta)=\mrm{U}$,
\item for every $n\in\NN,\xi_1,\dotsc,\xi_n\in \Dom(P^{-1/4})\cap \Dom(P^{1/4})$ the matrix $\bigl((\omega_{\xi_i,\xi_j}\tau_{i/4}\otimes\id)(\mrm{U}^*)\bigr)_{i,j=1}^{n}$ preserves the positive cone in the sense of Definition \ref{def1}.
\end{enumerate}
\end{theorem}

\begin{remark}\label{rem3}\noindent
\begin{itemize}
\item Since $\tau_t(x)=P^{it} x P^{-it}$ for $x\in \LL^{\infty}(\GG),t\in\RR$, the functional $\omega_{\xi_i,\xi_j}\tau_{i/4}$ in Theorem \ref{thm1} is well defined in $\LL^1(\GG)$ as $\omega_{P^{-1/4} \xi_i,P^{1/4}\xi_j}$.
\item The easier implication of Theorem \ref{thm1}, i.e.~that $\mrm{U}^{\alpha}$ satisfies $(1)$--$(3)$, was essentially established as \cite[Proposition 4.10]{DawsSkalskiViselter}. In the case of discrete quantum groups the converse was proven in \cite[Proposition 4.11 (c)]{DawsSkalskiViselter}. Our contribution is the converse statement, that conditions $(1)$--$(3)$ for $\mrm{U}$ imply $\mrm{U}=\mrm{U}^\alpha$ in locally compact case.
\item There are minor differences between condition $(3)$ and the condition in \cite{DawsSkalskiViselter}.
\begin{itemize}
\item Both operators $P^{it}$ and $\nabla_{\wh\vp}^{it}$ implement scaling automorphism $\tau_t$ \cite[Proposition 2.1]{KustermansVaes}. We have decided to use operator $P$ rather than $\nabla_{\widehat{\vp}}$ as the authors of \cite{DawsSkalskiViselter}, since unitary $P^{it}$ is the standard implementation of $\tau_t$ in the sense of \cite[Theorem 2.3]{HaagerupStandard}, when we treat $\LL^2(\GG)$ as the GNS representation of $\vp$ (or $\psi$) -- see Example \ref{ex6}.
\item  In \cite{DawsSkalskiViselter}, unitary $\mrm{U}$ rather than $\mrm{U}^*$ appears. This change is not substantial, as applying entry-wise adjoint gives equivalent condition -- the positive cone is self-dual, see \cite[Remark 4.9]{DawsSkalskiViselter}.
\item  We prefer to express the third condition using matrix $((\omega_{\xi_i,\xi_j}\tau_{i/4}\otimes\id)(\mrm{U}^*))_{i,j=1}^{n}$, as opposed to e.g.~$((\omega_{\nabla_{\wh\vp}^{-1/2}\xi_i,\xi_j}\otimes\id)(\mrm{U}^*))_{i,j=1}^{n}$ (for appropriate vectors). This condition is arguably more natural/symmetric, as $(\omega_{\xi_i,\xi_j})_{i,j=1}^{n}$ is CP in $\M_n(\LL^1(\GG))\subseteq \CB(\LL^{\infty}(\GG),\M_n)$. Moreover, this formulation admits obvious formulation at the universal level, see Theorem \ref{thm5} and Definition \ref{def2}. Both conditions are equivalent as we can approximate with smooth vectors for both $\nabla_{\wh\vp}$ and $P$ (these operators strongly commute), and then $\omega_{\xi_i,\xi_j}\tau_{i/4}=\omega_{\nabla_{\wh\vp}^{-1/2} (\nabla_{\wh\vp}^{1/4}\xi_i),\nabla_{\wh\vp}^{1/4}\xi_j}$.
\end{itemize}
\item More substantially, the authors of \cite{DawsSkalskiViselter} do not need to assume compatibility with modular conjugation. However, this condition is important for our proof.
\item Our reasoning is considerably different from the argument in \cite{DawsSkalskiViselter}.
\end{itemize}
\end{remark}

We will discuss several examples of the standard implementation after the proof.

\begin{proof}[Proof of Theorem \ref{thm1}]
If $\mrm{U}=\mrm{U}^\alpha$ is the standard implementation, then conditions $(1)$ and $(2)$ are clear and condition $(3)$ was proved in \cite[Proposition 4.10]{DawsSkalskiViselter} (see Remark \ref{rem3}). Consequently, assume that $\mrm{U}\in \LL^{\infty}(\GG)\bar\otimes \B(\msf{H}_\theta)$ is a unitary representation of $\GG$ satisfying $(1)$--$(3)$.\\

Let us introduce two operators on $\LL^2(\GG)\otimes\msf{H}_\theta$:
\begin{equation}\label{eq27}
J_{\mrm{U}}=(J_{\wh\vp}\otimes J_\theta)\mrm{U},\qquad
A_\mrm{U}=J_\mrm{U} J_{\wh\theta} \nabla_{\wh\theta}^{1/2},\quad\Dom(A_\mrm{U})=\Dom(\nabla_{\wh\theta}^{1/2}).
\end{equation}
Operator $J_\mrm{U}$ is anti-unitary, which is self-adjoint by the second assumption. Operator $A_\mrm{U}$ is densely defined, linear and closed. Our aim is to prove $J_\mrm{U}=J_{\wh\theta}$ and $A_\mrm{U}=\nabla_{\wh\theta}^{1/2}$. To do this, we will show that $A_\mrm{U}$ is positive, self-adjoint and use uniqueness of polar decomposition. We need to recall some useful facts concerning weight $\wh\theta\in \oon{W}_{\nsf}(\GG\ltimes\N)$. Its GNS map is quite concrete:
\begin{equation}\label{eq41}
\lin\{(\hat{x}\otimes \I)\alpha(y)\mid
\hat{x}\in \mf{N}_{\wh\vp},y\in \mf{N}_\theta\}\subseteq \mf{N}_{\wh\theta},\quad 
\Lambda_{\wh\theta}\bigl(
(\hat{x}\otimes \I)\alpha(y)\bigr)=
\Lambda_{\wh\vp}(\hat{x})\otimes \Lambda_\theta(y)
\end{equation}
and the above subspace is a $\sigma$-$\operatorname{SOT}^*\times\|\cdot\|$ core for $\Lambda_{\wh\theta}$ \cite[Page 439]{VaesUnitary}. Furthermore, 
\begin{equation}\label{eq28}\begin{split}
\lin \{\alpha(y^*)(\xi\otimes \Lambda_\theta(z))\mid \xi\in & \Dom(\nabla_{\wh\vp}^{1/2}),y,z\in\mf{N}_\theta\}
\subseteq \Dom(\nabla_{\wh\theta}^{1/2}),\\
J_{\wh\theta}\nabla_{\wh\theta}^{1/2}\bigl(
\alpha(y^*)(\xi\otimes \Lambda_\theta(z))
\bigr)&=
\alpha(z^*)(J_{\wh\vp}\nabla_{\wh\vp}^{1/2}\xi\otimes \Lambda_\theta(y))
\end{split}\end{equation}
and the above space is a (norm) core for $J_{\wh\theta}\nabla_{\wh\theta}^{1/2}$ \cite[Lemma 3.11]{VaesUnitary}. As $P$ and $\nabla_{\wh\vp}$ strongly commute, the smaller space 
\[
\lin \{\alpha(y^*)(\xi\otimes \Lambda_\theta(z))\mid \xi\in \bigcup_{\eps>0} 1_{\left[\eps,\eps^{-1}\right]}( \nabla_{\wh\vp})\,1_{\left[\eps,\eps^{-1}\right]}( P)\LL^2(\GG),y,z\in\mf{N}_\theta\}
\]
also forms a core for $J_{\wh\theta} \nabla_{\wh\theta}^{1/2}$. First we claim that
\begin{equation}\label{eq29}
\la \zeta | A_\mrm{U} \zeta \ra \ge 0
\end{equation}
for $\zeta=\sum_{k=1}^{K} \alpha(y_k^*)(\xi_k\otimes \Lambda_\theta(z_k))$, with $\xi_k\in \bigcup_{\eps>0} 1_{\left[\eps,\eps^{-1}\right]}( \nabla_{\wh\vp})\,1_{\left[\eps,\eps^{-1}\right]}( P)\LL^2(\GG), y_k,z_k\in \mf{N}_\theta$. Using \eqref{eq28}, \cite[Equation (5), Page 16]{Stratila} and the assumptions on $\mrm{U}$ we calculate
\[\begin{split}
&\quad\;
\la \zeta|A_\mrm{U}\zeta\ra =
\sum_{n,k=1}^{K} 
\bigl\la 
\alpha(y_n^*)(\xi_n\otimes \Lambda_\theta(z_n)) \big|
(J_{\wh\vp}\otimes J_\theta)\mrm{U}
J_{\wh\theta}\nabla_{\wh\theta}^{1/2}
\bigl(\alpha(y_k^*)(\xi_k\otimes \Lambda_\theta(z_k))\bigr)\bigr\ra \\
&=
\sum_{n,k=1}^{K} 
\bigl\la 
\alpha(y_n^*)(\xi_n\otimes \Lambda_\theta(z_n)) \big|
(J_{\wh\vp}\otimes J_\theta)\mrm{U}
\alpha(z_k^*)(J_{\wh\vp}\nabla_{\wh\vp}^{1/2}\xi_k\otimes \Lambda_\theta(y_k))\bigr\ra \\
&=
\sum_{n,k=1}^{K} 
\bigl\la 
\alpha(y_n^*)(\xi_n\otimes \Lambda_\theta(z_n)) \big|
(\I\otimes J_\theta z_k^* J_\theta)
(J_{\wh\vp}\otimes J_\theta)\mrm{U}(J_{\wh\vp}\nabla_{\wh\vp}^{1/2}\xi_k\otimes \Lambda_\theta(y_k))\bigr\ra \\
&=
\sum_{n,k=1}^{K} 
\bigl\la \xi_n\otimes J_\theta z_k J_\theta \Lambda_\theta(z_n)  \big|
\alpha(y_n)
(J_{\wh\vp}\otimes J_\theta)
\mrm{U}(J_{\wh\vp}\nabla_{\wh\vp}^{1/2}\xi_k\otimes \Lambda_\theta(y_k))\bigr\ra \\
&=
\sum_{n,k=1}^{K} 
\bigl\la \xi_n\otimes J_\theta z_k J_\theta \Lambda_\theta(z_n)  \big|
\mrm{U}^*(\I\otimes y_n ) \mrm{U}
\mrm{U}^*(\nabla_{\wh\vp}^{1/2}\xi_k\otimes J_\theta \Lambda_\theta(y_k))\bigr\ra \\
&=
\sum_{n,k=1}^{K} 
\bigl\la J_\theta z_k J_\theta \Lambda_\theta(z_n)  \big|
(\omega_{\xi_n, \nabla_{\wh\vp}^{1/2}\xi_k}\otimes\id)
(\mrm{U}^*)  y_n J_\theta \Lambda_\theta(y_k)\bigr\ra \\
&=
\sum_{n,k=1}^{K}\bigl\la 
z_n J_\theta \Lambda_\theta(z_k) \Big|
(\omega_{\nabla_{\wh\vp}^{-1/4} (\nabla_{\wh\vp}^{1/4} \xi_n),\nabla_{\wh\vp}^{1/4}( \nabla_{\wh\vp}^{1/4}\xi_k)}\otimes\id)
(\mrm{U}^*)  y_n J_\theta \Lambda_\theta(y_k)\bigr\ra\\
&=
\bigl\la \bigl(
z_n J_\theta \Lambda_\theta(z_k) \bigr)_{n,k=1}^{K} \Big|
\bigl(
(\omega_{\nabla_{\wh\vp}^{1/4}\xi_n ,\nabla_{\wh\vp}^{1/4} \xi_k}\tau_{i/4}\otimes\id)
(\mrm{U}^*)  y_n J_\theta \Lambda_\theta(y_k)\bigr)_{n,k=1}^{K}\bigr\ra
\end{split}\]
and Lemma \ref{lemma3} gives $\la \zeta|A_\mrm{U}\zeta\ra\ge 0$ (cone $\mf{P}_{\Tr_K\otimes \theta}$ is self-dual). Since the space of vectors $\zeta$ as above forms a core for $J_{\wh\theta}\nabla_{\wh\theta}^{1/2}$, it is also a core for $A_\mrm{U}$ and we conclude that $A_\mrm{U}$ is positive.

It follows that $A_\mrm{U}$ is symmetric, i.e.~$A_\mrm{U}\subseteq A_\mrm{U}^*$. Using \cite[Proposition 8.14]{SoltanPrimer} we see $A_\mrm{U}^*=\nabla_{\wh\theta}^{1/2} J_{\wh\theta}J_\mrm{U}$, hence
\begin{equation}\label{eq30}
\Dom(\nabla_{\wh\theta}^{1/2})=\Dom(A_\mrm{U})\subseteq \Dom(A_\mrm{U}^*)=
J_\mrm{U} J_{\wh\theta} \Dom(\nabla_{\wh\theta}^{1/2})=w \Dom(\nabla_{\wh\theta}^{1/2}),
\end{equation}
where we denote $w=J_\mrm{U} J_{\wh\theta}$. Observe that
\[
w=J_\mrm{U} J_{\wh\theta}=(J_{\wh\vp}\otimes J_\theta)\mrm{U} J_{\wh\theta}=
\mrm{U}^*(J_{\wh\vp}\otimes J_\theta)J_{\wh\theta}=\mrm{U}^*\mrm{U}^{\alpha} \in \LL^{\infty}(\GG)\bar\otimes \B(\msf{H}_\theta).
\]
 Moreover
\[\begin{split}
&\quad\;
w(\I\otimes J_\theta n J_\theta)w^*=
\mrm{U}^* \mrm{U}^\alpha(J_{\wh\vp}J_{\wh\vp}\otimes J_\theta n J_\theta)\mrm{U}^{\alpha *} \mrm{U}=
\mrm{U}^* (J_{\wh\vp}\otimes J_\theta) \mrm{U}^{\alpha *}(\I\otimes n)
\mrm{U}^\alpha (J_{\wh\vp}\otimes J_\theta)\mrm{U}\\
&=
(J_{\wh\vp}\otimes J_\theta) \mrm{U}\alpha(n)\mrm{U}^*
 (J_{\wh\vp}\otimes J_\theta)=\I\otimes J_\theta n J_\theta
\end{split}\]
for $n\in\N$, hence $w$ is a unitary in $ \LL^{\infty}(\GG)\bar\otimes \N$. It is an $\alpha$-cocycle:
\[
(\Delta\otimes\id)(w)=
\mrm{U}^*_{[23]} \mrm{U}^*_{[13]} \mrm{U}^{\alpha}_{[13]} \mrm{U}^\alpha_{[23]}=\mrm{U}^*_{[23]} \mrm{U}^*_{[13]} \mrm{U}^{\alpha}_{[13]} \mrm{U}_{[23]}\mrm{U}^*_{[23]}\mrm{U}^\alpha_{[23]}=
(\id\otimes \alpha)(w)w_{[23]}.
\]
One checks on generators that $\beta=\Ad(w)$ restricts to an automorphism of $\GG\ltimes\N$ with inverse $\beta^{-1}=\Ad(w^*)\in \Aut(\GG\ltimes\N)$. Moreover, 
\[
\beta(z),\;\beta^{-1}(z)\in \mf{N}_{\wh\theta},\quad 
\Lambda_{\wh\theta}(\beta(z))=w \Lambda_{\wh\theta}(z),\;
\Lambda_{\wh\theta}(\beta^{-1}(z))=w^* \Lambda_{\wh\theta}(z)\quad(z\in\mf{N}_{\wh\theta})
\]
by the proof of \cite[Proposition 4.2]{VaesUnitary}. Using this we obtain
\[
J_{\wh\theta}\nabla_{\wh\theta}^{1/2} w \Lambda_{\wh\theta}(z)=
J_{\wh\theta}\nabla_{\wh\theta}^{1/2} \Lambda_{\wh\theta}(\beta(z))=
\Lambda_{\wh\theta}(\beta(z^*))=
w \Lambda_{\wh\theta}(z^*)=
wJ_{\wh\theta}\nabla_{\wh\theta}^{1/2} \Lambda_{\wh\theta}(z)
\]
for $z\in \mf{N}_{\wh\theta}\cap {\mf{N}_{\wh\theta}}^*$, and as $\Lambda_{\wh\theta}(\mf{N}_{\wh\theta}\cap{\mf{N}_{\wh\theta}}^*)$ is a core for $J_{\wh\theta}\nabla_{\wh\theta}^{1/2}$ we deduce $w J_{\wh\theta} \nabla_{\wh\theta}^{1/2} \subseteq J_{\wh\theta} \nabla_{\wh\theta}^{1/2}w $, and $J_{\wh\theta} \nabla_{\wh\theta}^{1/2} \subseteq w^* J_{\wh\theta} \nabla_{\wh\theta}^{1/2}w $. Similar reasoning for $\beta^{-1}$ gives the converse inclusion, hence
\[
J_{\wh\theta} \nabla_{\wh\theta}^{1/2} =w^* J_{\wh\theta} \nabla_{\wh\theta}^{1/2}w =
w^* J_{\wh\theta} w\, w^* \nabla_{\wh\theta}^{1/2} w.
\]
Uniqueness of the polar decomposition implies
\[
\nabla_{\wh\theta}^{1/2} = w \nabla_{\wh\theta}^{1/2} w^*\quad \Rightarrow\quad 
\Dom(\nabla_{\wh\theta}^{1/2})=
w \Dom(\nabla_{\wh\theta}^{1/2})
\]
and \eqref{eq30} shows that $A_\mrm{U}$ is self-adjoint. Using again uniqueness of polar decomposition with equality \eqref{eq27} gives $J_\mrm{U} J_{\wh\theta}=\I$ and $J_{\wh\theta}=(J_{\wh\vp}\otimes J_\theta)\mrm{U}$, i.e.~$\mrm{U}=(J_{\wh\vp}\otimes J_\theta)J_{\wh\theta}=\mrm{U}^{\alpha}$.
\end{proof}

Let us now discuss some examples of the standard implementation.

\begin{example}\label{ex1}
While in general there is no ``formula'' for the standard implementation, we have a useful concrete expression if the weight $\theta$ is relatively $\rho$-invariant for a group-like $\rho$ affiliated with $\LL^{\infty}(\GG)$ (for an explanation of this terminology see Section \ref{sec:RIW} and in particular equation \eqref{RI1}). In this special case one has
\[
\mrm{U}^{\alpha *} \bigl(\rho^{1/2}\Lambda_{\vp}(x) \otimes \Lambda_\theta(n)\bigr)=
\Lambda_{\vp\otimes\theta}\bigl(\alpha(n)(x\otimes\I)\bigr)\quad(x\in \mf{N}_{\vp},\Lambda_\vp(x)\in \Dom(\rho^{1/2}),n\in\mf{N}_\theta)
\]
and 
\[
(\omega_{\xi,\eta}\otimes\id)(\mrm{U}^{\alpha *})  \Lambda_\theta(n)=
\Lambda_{ \theta}\bigl((\omega_{\xi,\rho^{-1/2}\eta}\otimes\id)\alpha(n) \bigr)\quad(
\xi\in \LL^2(\GG),\eta\in \Dom(\rho^{-1/2}),
n\in\mf{N}_\theta)
\]
(\cite[Proposition 4.3]{VaesUnitary}, \cite[Lemma 7.5]{DeCommerKrajczok}).
\end{example}

\begin{example}\label{ex2}
If $(\N,\alpha,\theta)$ is the von Neumann algebra $\LL^{\infty}(\GG)$ with translation action $\Delta$ and left Haar weight $\vp$, then $\mrm{U}^{\alpha}=\ww^{\GG}$, by Example \ref{ex1}. We obtain the same implementation if we take the right Haar integral instead. This is a consequence of Remark \ref{rem1} -- the unitary conjugation $u$ in this situation is an imaginary power of the scaling constant, hence it does not change the unitary implementation. 
\end{example}

\begin{example}\label{ex2.5}
Let $(\N,\alpha)$ be a von Neumann algebra with an action of $\GG$ and let $\theta\in \oon{W}_{\nsf}(\N)$ be a weight with dual weight $\wh{\theta}\in \oon{W}_{\nsf}(\GG\ltimes \N)$. As explained at the beginning of this section, one can identify the GNS representation with $\GG\ltimes \N\subseteq \B(\LL^2(\GG)\otimes\msf{H}_\theta)$. Crossed product carries a dual action $\alpha^{\ltimes}\colon \whG^{\op}\curvearrowright \GG\ltimes \N$, uniquely characterised by
\[
\alpha^{\ltimes}(\alpha(n))=\I\otimes \alpha(n),\quad \alpha^{\ltimes}(\hat{x}\otimes \I)=\Delta_{\whG^{\op}}(\hat{x})\otimes \I\quad( n\in\N,\hat{x}\in \LL^{\infty}(\whG))
\]
\cite[Proposition 2.2]{VaesUnitary}. Under this action, weight $\wh\theta$ is relatively $\delta_{\whG}$-invariant \cite[Proposition 3.10]{VaesUnitary} and the standard implementation is $\ww^{\whG^{\op}}_{[12]}$. Indeed, we have
\begin{equation}\label{eq31}\begin{split}
&\quad\;
\ww^{\whG^{\op} *}_{[12]}\bigl(\delta_{\whG}^{1/2} \Lambda_{\vp_{\whG^{\op}}}(\hat{y})\otimes \Lambda_{\wh\theta}((\hat{x}\otimes\I)\alpha(n))\bigr)=
\ww^{\whG^{\op} *}_{[12]}\bigl(\delta_{\whG}^{1/2} \Lambda_{\vp_{\whG^{\op}}}(\hat{y})\otimes \Lambda_{\vp_{\whG}}(\hat{x})\otimes \Lambda_{\theta}(n)\bigr)\\
&=
\ww^{\whG^{\op} *}_{[12]}\bigl(\delta_{\whG^{\op}}^{-1/2} \Lambda_{\vp_{\whG^{\op}}}(\hat{y})\otimes \Lambda_{\psi_{\whG^{\op}}}(\hat{x})\otimes \Lambda_{\theta}(n)\bigr)=
 \Lambda_{\vp_{\whG^{\op}}\otimes \psi_{\whG^{\op}}}\bigl( 
 \Delta_{\whG^{\op}}(\hat{x})(\hat{y}\otimes\I)\bigr)\otimes \Lambda_{\theta}(n)\\
 &
=
 \Lambda_{\vp_{\whG^{\op}}\otimes \vp_{\whG}}\bigl( 
 \Delta_{\whG^{\op}}(\hat{x})(\hat{y}\otimes\I)\bigr)\otimes \Lambda_{\theta}(n) \overset{\star}{=}
 \Lambda_{\vp_{\whG^{\op}}\otimes\wh\theta}\bigl(
\alpha^{\ltimes}\bigl((\hat{x}\otimes\I)\alpha(n)\bigr)\,
 (\hat{y}\otimes\I\otimes\I)
 \bigr)\\
 &=
 \mrm{U}^{\wh\alpha *}\bigl(\delta_{\whG}^{1/2} \Lambda_{\vp_{\whG^{\op}}}(\hat{y})\otimes \Lambda_{\wh\theta}((\hat{x}\otimes\I)\alpha(n))\bigr)
\end{split}\end{equation}
for $\hat{x}\in \mf{N}_{\vp_{\whG}}, \hat{y}\in \mf{N}_{\vp_{\whG^{\op}}},n\in\mf{N}_\theta$ such that $\Lambda_{\vp_{\whG^{\op}}}(\hat{y})\in\Dom(\delta_{\whG}^{1/2})$. In this calculation we have used \eqref{eq41}, $\vp_{\whG}=\psi_{\whG^{\op}}$, $\delta_{\whG^{\op}}=\delta_{\whG}^{-1}$, Example \ref{ex1} and that right Haar integral $\psi_{\HH}$ is relatively $\delta_{\HH}^{-1}$-invariant under translation action (see Example \ref{ex7}).  Equality $\overset{\star}{=}$ follows from
\[\begin{split}
&\quad\;(\la \Lambda_{\vp_{\whG^{\op}}}(\hat{z})|\otimes\id\otimes\id)\bigl(
 \Lambda_{\vp_{\whG^{\op}}\otimes \vp_{\whG}}\bigl( 
 \Delta_{\whG^{\op}}(\hat{x})(\hat{y}\otimes\I)\bigr)\otimes \Lambda_{\theta}(n) \bigr)\\
 &=
 \Lambda_{\vp_{\whG}}\bigl(
 (\omega_{\Lambda_{\vp_{\whG^{\op}}}(\hat{z}),\Lambda_{\vp_{\whG^{\op}}}(\hat{y})}\otimes\id)\Delta_{\whG^{\op}}(\hat{x})\bigr)\otimes \Lambda_\theta(n)\\
 &=
  \Lambda_{\wh\theta}\bigl((
 (\omega_{\Lambda_{\vp_{\whG^{\op}}}(\hat{z}),\Lambda_{\vp_{\whG^{\op}}}(\hat{y})}\otimes\id)\Delta_{\whG^{\op}}(\hat{x})\otimes \I)
\alpha(n)
 \bigr)\\
 &=
 (\la \Lambda_{\vp_{\whG^{\op}}}(\hat{z})|\otimes\id\otimes\id)\Lambda_{\vp_{\whG^{\op}}\otimes\wh\theta}\bigl(
(\Delta_{\whG^{\op}}(\hat{x})\otimes\I)
(\I\otimes \alpha(n))
(\hat{y}\otimes\I\otimes\I)
 \bigr)\\
 &=
 (\la \Lambda_{\vp_{\whG^{\op}}}(\hat{z})|\otimes\id\otimes\id)\Lambda_{\vp_{\whG^{\op}}\otimes\wh\theta}\bigl(
\alpha^{\ltimes}\bigl((\hat{x}\otimes\I)\alpha(n)\bigr)(\hat{y}\otimes\I\otimes\I)
 \bigr)
\end{split}\]
for $\hat{z}\in \mf{N}_{\psi_{\whG}}$. Relation $\ww^{\whG^{\op}}_{[12]}=\mrm{U}^{\wh\alpha}$ follows now from \eqref{eq31} and density argument. Compare also \cite[Theorem 1.3]{DeCommerDeRo}.
\end{example}

\begin{example}\label{ex3}
Assume that $(\GG,\wh{\oon{R}})$ is quasi-triangular LCQG acting on von Neumann algebras $\alpha^{\M}\colon \GG\curvearrowright\M,  \alpha^{\N}\colon \GG\curvearrowright\N$. Then we can build braided tensor product $\M\ov\boxtimes \N$ with diagonal action $\alpha^{\M\ov\boxtimes\N}\colon \GG\curvearrowright \M\ov\boxtimes \N$ (for the relevant definitions and construction, see \cite{DeCommerKrajczok}). \cite[Theorem 5.11]{DeCommerKrajczokStandard} shows that the defining representation $\M\ov\boxtimes \N\subseteq \B(\LL^2(\M)\otimes\LL^2(\N))$ is standard\footnote{ With particular modular conjugation and unspecified positive cone; recall that standard implementation depends only on this data -- not the weight itself (Remark \ref{rem1}).} and the standard implementation of $\alpha^{\M\ov\boxtimes\N}$ is the tensor product of relevant implementations $\mrm{U}^{\M},\mrm{U}^{\N}$.
\end{example}

Next we have a look at the amplified action.

\begin{example}\label{ex4}
Let $\OnP$ be a von Neumann algebra and $\alpha\colon \GG\curvearrowright \N$ be an action of a (general) LCQG on a von Neumann algebra. Choose weights $\theta_{\OnP}\in \oon{W}_{\nsf}(\OnP), \theta_{\N}\in \oon{W}_{\nsf}(\N)$ and let $\mrm{U}^{\alpha}$ be the corresponding standard implementation of $\alpha$. Then $\beta= \alpha\otimes\id$ is an action $\GG\curvearrowright \N\bar\otimes \OnP$. We claim that the standard implementation of $\beta$ constructed with respect to $\theta_{\N}\otimes \theta_{\OnP}$ is $\mrm{U}^\beta=\mrm{U}^{\alpha}_{[12]}$.

We can identify $\msf{H}_{\theta_{\N}\otimes  \theta_{\OnP}}=\msf{H}_{\theta_{\N}}\otimes \msf{H}_{\theta_{\OnP}}$, and as the modular conjugation is given by $J_{\theta_{\N}}\otimes J_{\theta_{\OnP}}$, we easily check that $\mrm{U}^{\alpha}_{[12]}$ satisfies conditions $(1)$ and $(2)$ of Theorem \ref{thm1}. To see the third condition, take $n\in\NN, \xi_1,\dotsc,\xi_n\in \Dom(P^{-1/4})\cap\Dom(P^{1/4})$ and $\zeta\in \mf{P}_{\Tr_n \otimes \theta_{\N}\otimes \theta_{\OnP}}$. Arguing by approximation, we can use \cite[Corollary 2.9]{MiuraTomiyama} and assume that $\zeta=\sum_{a,a'=1}^{A} \zeta_{ 1,a,a'}\otimes \zeta_{2,a,a'}$, where matrix $(\zeta_{1,a,a'})_{a,a'=1}^{A}$ belongs to $\mf{P}_{\Tr_A\otimes \Tr_n\otimes \theta_{\N}}=\mf{P}_{\Tr_{A\times n}\otimes \theta_{\N}}$ and $(\zeta_{2,a,a'})_{a,a'=1}^{A}$ belongs to $\mf{P}_{\Tr_A\otimes \theta_{\OnP}}$ (\cite[Definition 2.2, Proposition 2.4]{MiuraTomiyama}). Since $\mrm{U}^{\alpha}$ is standard, matrix $((\omega_{\xi_{a,i},\xi_{a',i'}}\tau_{i/4}\otimes \id)(\mrm{U}^{\alpha *}))_{(a,i),(a',i')=(1,1)}^{A\times n}$ with $\xi_{a,i}=\xi_i$ preserves the positive cone $\mf{P}_{\Tr_{A\times n}\otimes \theta_{\N}}$. The claim follows now from
\[
\bigl((\omega_{\xi_i,\xi_j}\tau_{i/4}\otimes\id)(\mrm{U}^{\alpha *}_{[12]})\bigr)_{i,j=1}^{n} \bullet \zeta =
\sum_{a,a'=1}^{A} \wt{\zeta}_{1,a,a'}\otimes \zeta_{2,a,a'}
\]
(writing $\bullet$ for the Schur multiplication), where
\[
(\wt{\zeta}_{1,a,a'})_{a,a'=1}^A=
\bigl((\omega_{\xi_{a,i},\xi_{a',i'}}\tau_{i/4}\otimes\id)(\mrm{U}^{\alpha *})\bigr)_{(a,i),(a',i')=(1,1)}^{A\times n}\bullet 
(\zeta_{1,a,a'})_{a,a'=1}^A\in \mf{P}_{\Tr_{A\times n}\otimes \theta_{\N}}
\]
Passing via the flip isomorphisms, we obtain a similar conclusion for $\GG\curvearrowright \OnP\bar\otimes \N$.
\end{example}

The next example we discuss is the case of classical locally compact groups. Haagerup in \cite[Corollary 3.6]{HaagerupStandard} showed that their actions admit a canonical implementation, we check that it corresponds to the standard implementation, as understood in this article.

\begin{example}\label{ex5}
Let $G$ be a locally compact group with a (point-weak$^*$-continuous) action on a von Neumann algebra $\N$ equipped with $\theta\in \oon{W}_{\nsf}(\N)$. For $g\in G$, let $u_g\in \B(\msf{H}_\theta)$ be the standard implementation of $\alpha_g\in \Aut(\N)$ in the sense of 
\cite[Theorem 3.2]{HaagerupStandard}, in particular $\alpha_g(x)=u_g x u_g^*$ for $g\in G,x\in \N$. Using the identification $\LL^2(G)\otimes \msf{H}_\theta\simeq \LL^2(G;\msf{H}_\theta)$ (cf.~\cite[Proposition E.2]{Cohn}) one can check that the function $g\mapsto u_g$ defines an element $\mrm{U}\in \LL^{\infty}(G)\bar\otimes \B(\msf{H}_\theta)$ and $x\mapsto (g\mapsto \alpha_{g^{-1}}(x))$ a left action $G\curvearrowright \N$. With Lemma \ref{lemma3}, it is easy to see that $\mrm{U}$ satisfies conditions from Theorem \ref{thm1}, hence $\mrm{U}=\mrm{U}^\alpha$.
\end{example}

\begin{example}\label{ex6}
As a special case of Example \ref{ex5}, consider locally compact quantum group $\GG$ and action $\alpha\colon \RR\curvearrowright \LL^{\infty}(\GG)$ given by the scaling group (i.e.~$\alpha(x)=(t\mapsto \tau_{-t}(x))$; we have the minus sign as we work with left actions) or modular automorphism group of $\vp$ or $\psi$. Consider $\LL^{\infty}(\GG)\subseteq \B(\LL^2(\GG))$ as the GNS representation of $\vp$ or $\psi$. In all these cases, the standard implementation of $\RR\curvearrowright \LL^{\infty}(\GG)$ is given by the family of operators $P^{it},\nabla_\vp^{it}$ or $\nabla_\psi^{it}\,(t\in\RR)$. This is again a consequence of Theorem \ref{thm1} and the observation that operators $P^{it},\nabla_\vp^{it}, \nabla_\psi^{it}$ commute with $J_{\vp},J_{\psi}$ and preserve positive cones $\mf{P}_{\vp},\mf{P}_{\psi}$; see \cite[Proposition 2.13 and Page 74]{KustermansVaes}.
\end{example}

See also \cite[Theorem 5.1]{KasprzakKhosraviSoltan}, for an example of the standard implementation of an action associated to a morphism between quantum groups.

\section{Relatively invariant weights}\label{sec:RIW}

Let $\alpha\colon \GG\curvearrowright \N$ be an action of a locally compact quantum group $\GG$ on a von Neumann algebra $\N$. Consider a n.s.f.~weight $\theta\in \oon{W}_{\nsf}(\N)$ and a strictly positive, self-adjoint operator $\rho\,\eta\LL^{\infty}(\GG)$ which is group-like in the sense $\Delta(\rho)=\rho\otimes \rho$. We say that $\theta$ is \emph{relatively $\rho$-invariant} if
\begin{equation}\label{RI1}\tag{RI1}
\theta\bigl((\omega_\xi\otimes\id)\alpha(n)\bigr)=\|\rho^{1/2} \xi\|^2 \theta(n)\quad (\xi\in \Dom(\rho^{1/2}), n\in\mf{M}_\theta^+).
\end{equation}
This definition was introduced in \cite[Definition 2.3]{VaesUnitary}, note however that we restrict our attention to the situation where $\rho$ is group-like. Equality \eqref{RI1} is not the only possible condition expressing ``$\rho$-invariance'' of $\theta$ with respect to the action $\alpha$. We consider several stronger conditions: 
\begin{equation}\label{RI2}\tag{RI2}
\theta\bigl((\omega_\xi\otimes\id)\alpha(n)\bigr)=\|\rho^{1/2} \xi\|^2 \theta(n)\quad (\xi\in \Dom(\rho^{1/2}), n\in\N^+)
\end{equation}
(i.e.~in \eqref{RI2} we do not assume that $n$ is integrable) and
\begin{equation}\label{RI3}\tag{RI3}
(\id\otimes\id\otimes \theta)
(\id\otimes\alpha)(X)
=
(\id\otimes \theta)(X)\otimes \rho\quad(X\in (\OnP\bar\otimes\N)^+)
\end{equation}
as an equality in $(\OnP\bar\otimes \LL^{\infty}(\GG))^{+,\ext}$ for every von Neumann algebra $\OnP$ and finally
\begin{equation}\label{RI4}\tag{RI4}
(\id\otimes\id\otimes \theta)
(\id\otimes\alpha)(X)
=
(\id\otimes \theta)(X)\otimes \rho\quad(X\in (\OnP\bar\otimes\N)^{+,\ext})
\end{equation}
in $(\OnP\bar\otimes\LL^{\infty}(\GG))^{+,\ext}$ for all von Neumann algebras $\OnP$. The main result of this section is the following.

\begin{theorem}\label{thm2}
Conditions \eqref{RI1}, \eqref{RI2}, \eqref{RI3} and \eqref{RI4} are equivalent.
\end{theorem}

\begin{remark}\label{rem2}\noindent
\begin{itemize}
\item Kustermans and Vaes proved in \cite[Proposition 3.1]{KustermansVaes} that the left Haar integral $\vp\in \oon{W}_{\nsf}(\LL^{\infty}(\GG))$ satisfies condition \eqref{RI3} with respect to the action $\Delta$, while definition of left Haar integral states \eqref{RI1}. Their proof is very specific to the situation of Haar integral and translation action, and does not seem to generalise to arbitrary actions.
\item Stronger form of invariance \eqref{RI3}  for Haar integral is used crucially in the construction of the canonical operator valued weight $\N\rightarrow \N^{\alpha}$ \cite[Proposition 1.3]{VaesUnitary}. This operator valued weight is used to construct the dual weight on $\GG\ltimes \N$ \cite[Definition 3.1]{VaesUnitary}.
\item \cite{BrannanChirvasituViselter} call weights satisfying \eqref{RI2} with $\rho=\I$, ``completely invariant''. In the situation when $\HH\le \GG$ is a closed quantum subgroup, the authors of \cite{BrannanChirvasituViselter} prove that there is a completely invariant n.s.f.~weight on $\LL^{\infty}(\GG/\HH)$ if, and only if the modular element $\delta_{\GG}$ restricts to $\delta_{\HH}$ \cite[Theorem 4.4]{BrannanChirvasituViselter}.
\end{itemize}
\end{remark}

Before we present a proof of Theorem \ref{thm2}, let us discuss some examples of relatively invariant weights.

\begin{example}\label{ex8}
By definition \cite[Definition 1.1]{KustermansVaes}, the left Haar integral $\vp\in \oon{W}_{\nsf}(\LL^{\infty}(\GG))$ is relatively $\I$-invariant (satisfies \eqref{RI1} for $\rho=\I$), with respect to the translation action $\Delta\colon \GG\curvearrowright \LL^{\infty}(\GG)$.
\end{example}

\begin{example}\label{ex7}
The right Haar integral $\psi\in \oon{W}_{\nsf}(\LL^{\infty}(\GG))$ is relatively $\delta^{-1}$-invariant for $\Delta\colon\GG\curvearrowright\LL^{\infty}(\GG)$. One way of seeing it is realising $\psi$ as a dual weight, see the next example.
\end{example}

\begin{example}
 Let $\N$ be a von Neumann algebra with action $\alpha\colon \GG\curvearrowright \N$ and weight $\theta\in \oon{W}_{\nsf}(\N)$. The dual weight $\hat{\theta}\in \oon{W}_{\nsf}(\GG\ltimes \N)$ is relatively $\delta_{\whG}$-invariant with respect to the dual action $\whG^{\oon{op}}\curvearrowright \GG\ltimes \N$ (see the proof of \cite[Proposition 3.10]{VaesUnitary}).
\end{example}

\begin{example}
Garcia Guinto-Hall-Nelson introduced in \cite{GuintoHallNelson} the \emph{middle Haar weight} $\text{\tengmag{0.7}\Tando}=\vp_{\delta^{1/2}}\in \oon{W}_{\nsf}(\LL^{\infty}(\GG))$ and showed that it is relatively $\delta^{-1/2}$-invariant under the translation action $\Delta\colon \GG\curvearrowright \LL^{\infty}(\GG)$ \cite[Theorem 2.3]{GuintoHallNelson}. It is connected to almost unimodularity: by \cite[Theorem 4.2]{GuintoHallNelson} $\text{\tengmag{0.7}\Tando}$ is almost periodic if, and only if $\whG$ is almost unimodular.
\end{example}

The last example will use terminology and results from Section \ref{sec:univ}.

\begin{example}\label{ex9}
Let $\GG\rightarrow \HH$ be a morphism of quantum groups described by strong quantum homomorphism $\Pi\colon \mrm{C}_0^u(\HH)\rightarrow \M(\mrm{C}_0^u(\GG))$, $\alpha_{\HH}\colon \HH\curvearrowright \N$ an action and $\alpha_{\GG}\colon \GG\curvearrowright \N$ the corresponding action given by composition (see the paragraph before Corollary \ref{cor1}). Assume that $\theta\in \oon{W}_{\nsf}(\N)$ is relatively $\rho$-invariant under $\alpha_{\HH}$, where $\rho \eta \LL^{\infty}(\HH)$ is a strictly positive, self-adjoint group-like operator. Denote by $\rho^u$ its lift, affiliated with $\mrm{C}_0^u(\HH)$ (see Preliminaries). In this situation $\theta$ is relatively $\pi_{\GG}\Pi(\rho^u)$-invariant under $\alpha_{\GG}$. Indeed, using Proposition \ref{prop2}
\[
\theta\bigl( (\omega_\xi\otimes\id) \alpha_{\GG}(n)\bigr)=
\theta\bigl( (\omega_\xi \pi_{\GG} \Pi \otimes\id)\alpha^{u}_{\HH}(n)\bigr)=
\|\pi_{\GG}\Pi(\rho^u)^{1/2}\xi\|^2 \theta(n)
\]
for $n\in \mf{M}_\theta^+, \xi\in \Dom(\pi_{\GG}\Pi(\rho^u)^{1/2})$, which proves \eqref{RI1} for $\theta$ with respect to $\alpha_{\GG}$.
\end{example}

In order to prove Theorem \ref{thm2}, we need to establish three lemmas. Lemma \ref{lemma4} is a general useful result concerning relatively invariant weight; we will use it in the proof of Theorem \ref{thm2} and in next sections. 

\begin{lemma}\label{lemma4}
Let $\alpha\colon \GG\curvearrowright \N$ be an action, $\rho\,\eta\LL^{\infty}(\GG)$ a strictly positive, self-adjoint group-like operator and $\theta\in \oon{W}_{\nsf}(\N)$ a relatively $\rho$-invariant weight.
\begin{itemize}
\item Operators $P$ and $\rho$ strongly commute.
\item  Write $\kappa_t=\Ad(P^{it} \rho^{-it})=\tau_t \Ad(\rho^{-it})\in \Aut(\LL^{\infty}(\GG))$ for $t\in\RR$. $\kappa$ is a continuous one parameter group of automorphisms.
\item We have $\alpha \sigma^\theta_t=(\kappa_t\otimes \sigma^\theta_t)\alpha$ for $t\in\RR$.
\item Write $\tau_t^{u,\vN}\in \Aut(\mrm{C}_0^u(\GG)^{**})$ for the normal extension of $\tau_t^u$ and define automorphism $\kappa_t^u=\tau^{u,\vN}_t\Ad(\rho^{u -it})\in \Aut (\mrm{C}_0^u(\GG)^{**})$. Then $\alpha^u\sigma^\theta_t=(\kappa_t^u\otimes \sigma^\theta_t)\alpha^u$ for $t\in\RR$.
\end{itemize}
\end{lemma}

\begin{proof}
These results are stated as a part of \cite[Proposition 2.4]{VaesUnitary} when $\rho$ is the inverse of modular element, with $P$ replaced by $\nabla_{\wh\vp}$ (up to differences in conventions). We can replace $\nabla_{\wh\vp}$ by $P$ as both operators implement scaling group on $\LL^{\infty}(\GG)$. The case of general $\rho$ is similar. See also \cite[Proposition 8.4, Corollary 8.5]{Kus02}. To see the last property, observe that
\[\begin{split}
\Delta^{u,r,\vN}\kappa_t(x)&=
\Delta^{u,r}( \rho^{-it}\tau_t(x)\rho^{it})=
(\rho^{u -it}\otimes \rho^{-it})(\tau^{u}_t\otimes\tau_t)\Delta^{u,r}(x)
(\rho^{u it}\otimes \rho^{it})\\
&=
\bigl(
\Ad (\rho^{u -it})\tau_t^{u,\vN} \otimes \Ad(\rho^{-it})\tau_t\bigr)
\Delta^{u,r,\vN}(x)
\end{split}\]
for $x\in \mrm{C}_0(\GG), t\in\RR$ thus $\Delta^{u,r,\vN}\kappa_t=(\kappa_t^u\otimes\kappa_t)\Delta^{u,r,\vN}$ by continuity and
\[\begin{split}
&\quad\;
(\id\otimes\alpha)\alpha^u\sigma^\theta_t
=
(\Delta^{u,r,\vN}\otimes \id)\alpha\sigma^\theta_t=
(\Delta^{u,r,\vN}\otimes \id)(\kappa_t\otimes \sigma^\theta_t)\alpha\\
&=
(\kappa_t^u\otimes\kappa_t\otimes \sigma^\theta_t)(\Delta^{u,r,\vN}\otimes \id)\alpha=
(\kappa_t^u\otimes\kappa_t\otimes \sigma^\theta_t)(\id\otimes\alpha)\alpha^u=
(\kappa_t^u\otimes \alpha \sigma^\theta_t)\alpha^u,
\end{split}\]
hence injectivity of $\id\otimes\alpha$ ends the proof.
\end{proof}

Next we establish two lemmas aimed at the proof of Theorem \ref{thm2}. Lemma \ref{lemma2} can be seen as an extension of ergodicity of $\Delta$ or ``invariants are constant'' technique \cite[Theorem 2.1]{DawsCPMults} to unbounded elements, incorporating also group-like $\rho$. Its proof is obtained by reduction to Lemma \ref{lemma1}, which is concerned with a simplified situation: trivial $\rho$, but more general action.

\begin{lemma}\label{lemma1}
Let $\GG$ be a locally compact quantum group, $\OnP\subseteq \B(\msf{H})$ a von Neumann algebra, $\mrm{U}\in \B(\msf{H})\bar\otimes \LL^{\infty}(\GG)$ a right representation of $\GG$ and $Y\in \OnP^{+,\ext}$. If 
\begin{equation}\label{eq8}
q_{Y\otimes \I}\bigl(\mrm{U}^*(\xi\otimes\eta)\bigr)=
q_Y(\xi)\|\eta\|^2
\end{equation}
for $\xi\in \Dom(q_Y), \eta\in \LL^2(\GG)$, then $\mrm{U}^*(Y\otimes \I) \mrm{U}=Y\otimes \I$.
\end{lemma}

\begin{proof}
Write $Y=A_Y+\infty (\I-p_Y)$ (see Preliminaries), then $Y\otimes \I=A_Y\otimes \I+\infty (\I\otimes \I-p_Y\otimes \I)$. Define Hilbert space $\mc{K}$ equal to $\Dom(Y^{1/2})=\Dom(A_Y^{1/2})$ with inner product
\[
\la \xi|\xi'\ra_{\mc{K}}=\la\xi|\xi'\ra_{\msf{H}}+
\la A_Y^{1/2} \xi|A_Y^{1/2}\xi'\ra_{\msf{H}}\quad(\xi,\xi'\in\mc{K}=\Dom(Y^{1/2})),
\]
let $\iota\colon \mc{K}\hookrightarrow \msf{H}$ be the canonical contractive inclusion, see \cite[Section VI.1.3]{Kato}. It is easy to check that for $\eta\in \mc{K}\otimes \LL^2(\GG)$ we have $(\iota\otimes\I)\eta\in \Dom(A_Y^{1/2}\otimes \I)$ and
\begin{equation}\label{eq9}
\|\eta\|^2_{\mc{K}\otimes\LL^2(\GG)}=
\|(\iota\otimes \I)\eta\|^2_{\msf{H}\otimes \LL^2(\GG)}+
\|(A_Y^{1/2}\otimes \I)(\iota\otimes\I)\eta\|^2_{\msf{H}\otimes\LL^2(\GG)}.
\end{equation}
Note that the above equation implies that $\iota\otimes \I$ is injective. Observe also
\begin{equation}\label{eq10}
\oon{Im}(\iota\otimes \I)=\Dom(A_Y^{1/2}\otimes \I)
\end{equation}
using \eqref{eq9}. For $\xi\in \mc{K}\odot \LL^2(\GG)$, define $w(\xi) \in\mc{K}\otimes \LL^2(\GG)$ via
\begin{equation}\label{eq11}
(\iota\otimes \I)w(\xi)=U^* (\iota\otimes \I)(\xi).
\end{equation}
This is well defined: using \eqref{eq8} (and extending these quadratic forms to sesquilinear forms) we have $q_{Y\otimes \I} \bigl(\mrm{U}^* (\iota\otimes\I)(\xi)\bigr)<+\infty$, hence $\mrm{U}^*(\iota\otimes\I)(\xi)\in \Dom(A_Y^{1/2}\otimes\I)=\oon{Im}(\iota\otimes\I)$ by \eqref{eq10}. Consequently, have a well defined linear map $w\colon \mc{K}\odot \LL^2(\GG)\rightarrow \mc{K}\otimes \LL^2(\GG)$. Since
\[\begin{split}
&\quad\;
\|w(\xi)\|^2_{\mc{K}\otimes \LL^2(\GG)}=
\|(\iota\otimes\I)w(\xi)\|^2+
\|(A_Y^{1/2}\otimes\I)(\iota\otimes\I)w(\xi)\|^2\\
&=
\|\mrm{U}^*(\iota\otimes\I)(\xi)\|^2+
\|(A_Y^{1/2}\otimes\I)\mrm{U}^*(\iota\otimes\I)(\xi)\|^2=
\|(\iota\otimes\I)\xi\|^2+
q_{Y\otimes \I}\bigl(\mrm{U}^*(\iota\otimes\I)(\xi)\bigr)\\
&=
\|(\iota\otimes\I)\xi\|^2+
q_{Y\otimes \I}\bigl((\iota\otimes\I)(\xi)\bigr)=\|\xi\|^2_{\mc{K}\otimes \LL^2(\GG)}
\end{split}\]
(where norms without subscript are calculated in $\msf{H}\otimes\LL^2(\GG)$) using \eqref{eq9}, \eqref{eq11}, \eqref{eq8}, $w$ extends to an isometry on $\mc{K}\otimes\LL^2(\GG)$. From \eqref{eq11} easily follows that $w\in \B(\mc{K})\bar\otimes \LL^\infty(\GG)$, and that $w$ is a right (isometric) representation of $\GG^{\oon{op}}$, i.e.~$\LL^{\infty}(\GG)$ with coproduct $\Delta(\cdot)_{[21]}$. By \cite[Corollary 4.16]{BrannanDawsSamei} $w$ is unitary. Equality \eqref{eq11} implies
\begin{equation}\label{eq17}
\mrm{U}^*(\iota\iota^*\otimes\I)\mrm{U}=
(\iota\otimes\I)w w^*(\iota^*\otimes \I)=
\iota\iota^*\otimes \I.
\end{equation}
Observe that
\begin{equation}\label{eq13}
\iota\iota^*(\eta)=0\;\textnormal{for}\;\eta\in \oon{Im}(p_Y)^{\perp}\quad\textnormal{ and}\quad\iota\iota^*(\eta)=(\I+A_Y)^{-1}\eta\;\textnormal{for}\; \eta\in \oon{Im}(p_Y).
\end{equation}
 Indeed, for $\xi\in \msf{H}$ we have
\[
\la \xi|\iota \iota^* \eta\ra=
\la \iota\iota^*  \xi|\eta\ra=
\la p_Y \iota\iota^*  \xi|\eta\ra=0
\]
in the first case and 
\[\begin{split}
&\quad\;
\la \xi | \iota\iota^*\eta\ra =
\la \iota \iota^* \xi |  \eta\ra=
\la \iota \iota^* \xi |  (\I+A_Y)^{-1}\eta\ra+
\la A_Y^{1/2}\iota \iota^* \xi |  A_Y^{1/2}(\I+A_Y)^{-1}\eta\ra\\
&=
\bigl\la \iota^*\xi|\iota^{-1}\bigl((\I+A_Y)^{-1}\eta\bigr)\bigr\ra_{\mc{K}}=
\bigl\la \xi|(\I+A_Y)^{-1}\eta\ra
\end{split}\]
in the second case, where we simply note that $(\I+A_Y)^{-1}\eta$ belongs to the image of $\iota$. This shows \eqref{eq13}. It follows in particular that the projection onto the kernel of $\iota\iota^*$ equals $1_{\{0\}}(\iota\iota^*)=\I-p_Y$ and equation \eqref{eq17} implies
\begin{equation}\label{eq19}
\mrm{U}^*(p_Y\otimes \I)\mrm{U}=p_Y\otimes \I,
\end{equation}
so that $\mrm{U}$ restricts to a unitary on $\oon{Im}(p_Y\otimes \I)$. Moreover $\iota\iota^*\rest_{\oon{Im} p_Y}=(\I+A_Y)^{-1}$, hence by \eqref{eq17}
\begin{equation}\label{eq43}
\mrm{U}^*( (\I+A_Y)^{-1}\otimes\I)\mrm{U}=(\I+A_Y)^{-1}\otimes \I\quad\textnormal{on}\quad \oon{Im}(p_Y\otimes \I).
\end{equation}
Now we claim that
\begin{equation}\label{eq18}
\mrm{U}^* ((Y+\I)\otimes \I)\mrm{U}=(Y+\I)\otimes \I\quad \textnormal{ in } \B(\msf{H}\otimes \LL^2(\GG))^{+,\ext}.
\end{equation}
Take $\xi\in \msf{H}\otimes \LL^2(\GG)$. If $((\I-p_Y)\otimes \I)\xi\neq 0$, then
\[
\omega_\xi((Y+\I)\otimes\I)\ge \infty \,\omega_\xi((\I-p_Y)\otimes \I)=\infty
\]
and similarly
\[
\omega_\xi\bigl(\mrm{U}^*((Y+\I)\otimes\I)\mrm{U}\bigr)\ge 
\infty\, \omega_\xi \bigl(\mrm{U}^* ((\I-p_Y)\otimes\I)\mrm{U}\bigr)=
\infty\, \omega_\xi \bigl((\I-p_Y)\otimes\I\bigr)=\infty
\]
using \eqref{eq19}. If $\xi \in \oon{Im}(p_Y\otimes\I)$, then
\[\begin{split}
&\quad\;
\omega_\xi((Y+\I)\otimes\I)=
q_{(A_Y+\I)\otimes\I}( (p_Y\otimes\I)\xi)=
q_{ (A_Y+\I)\otimes \I}(\mrm{U} (p_Y\otimes\I)\xi)\\
&=
q_{(A_Y+\I)\otimes\I}((p_Y\otimes\I)\mrm{U}\xi)=
\omega_{\mrm{U}\xi}((Y+\I)\otimes \I)=
\omega_{\xi}\bigl(\mrm{U}^*((Y+\I)\otimes \I)\mrm{U}\bigr)
\end{split}\]
using \eqref{eq43} in the second and \eqref{eq19} in the third equality. As arbitrary element of $\B(\msf{H}\otimes \LL^2(\GG))_*^+$ can be written as a series of vector functionals, this shows \eqref{eq18}, from which the claim follows.
\end{proof}

Next we upgrade Lemma \ref{lemma1} to incorporate group-like operator $\rho$, we also strengthen the conclusion.

\begin{lemma}\label{lemma2}
Let $\GG$ be a locally compact quantum group and $\rho\,\eta\,\LL^{\infty}(\GG)$ a strictly positive, self-adjoint operator satisfying $\Delta(\rho)=\rho\otimes \rho$. Let $\OnP\subseteq \B(\msf{H})$ be a von Neumann algebra and $Y\in (\OnP\bar\otimes\LL^{\infty}(\GG))^{+,\ext}$. If 
\begin{equation}\label{eq21}
q_{(\id\otimes \Delta)(Y)}(\xi\otimes\eta)=
q_Y(\xi)\|\rho^{1/2}\eta\|^2
\end{equation}
for $\xi\in \Dom(q_Y), \eta\in \Dom(\rho^{1/2})$, then there is a unique $C\in \OnP^{+,\ext}$ such that $Y=C\otimes \rho$.
\end{lemma}

\begin{proof}
First we rewrite the assumption to a more useful form:
\begin{equation}\label{eq22}
q_{Y\otimes \I}\bigl(\vv^*_{[23]}(\xi\otimes\eta)\bigr)=q_Y(\xi)\|\rho^{1/2}\eta\|^2\quad(\xi\in \Dom(q_Y), \eta\in\Dom(\rho^{1/2})).
\end{equation}
Indeed, write $Y=A_Y+\infty(\I-p_Y)$. Then for $\xi\in \Dom(q_Y), \eta\in \Dom(\rho^{1/2})$ by \eqref{eq21} we have $\xi\otimes\eta\in \oon{Im}((\id\otimes\Delta)(p_Y))$, hence $\vv_{[23]}^*(\xi\otimes\eta)\in \oon{Im}( p_Y\otimes\I)$ and
\[
q_{Y\otimes \I}\bigl(\vv^*_{[23]} (\xi\otimes\eta)\bigr)=
\|(A_Y^{1/2}\otimes\I)\vv^*_{[23]}(\xi\otimes\eta)\|^2=
\|(\id\otimes\Delta)(A_Y)^{1/2}(\xi\otimes\eta)\|^2=
q_{(\id\otimes\Delta)(Y)}(\xi\otimes\eta).
\]
Another application of \eqref{eq21} proves \eqref{eq22}.

The main part of the proof is a reduction to Lemma \ref{lemma1} by twisted amplification; this idea is similar in spirit to the Maharam extension in ergodic theory (see e.g.~\cite[Section 3.4]{Aaronson}). Choose two strictly positive, self-adjoint operators $\mathbf{A},\mathbf{B}$ on $\LL^2(\RR)$ satisfying $\mathbf{A}^{it}\mathbf{B}^{is}=e^{ist} \mathbf{B}^{is}\mathbf{A}^{it}\,(s,t\in\RR)$ (see e.g.~\cite[Section 2]{WoronowiczQEF}). Next, define unitary
\[
Z=\exp\bigl(-i\log(\mathbf{B})\otimes \log(\rho)\bigr)\in \B(\LL^2(\RR))\bar\otimes \LL^{\infty}(\GG)
\]
and observe that 
\begin{equation}\label{eq14}
Z^*(\mathbf{A}\otimes \rho)Z=\mathbf{A}\otimes \I.
\end{equation}
Indeed, \eqref{eq14} follows from
\[\begin{split}
&\quad\;
\bigl(Z^*(\mathbf{A}\otimes \rho)Z\bigr)^{it}=
Z^* (\mathbf{A}^{it}\otimes \rho^{it})Z=
Z^* \exp\bigl( -i \mathbf{A}^{it} \log(\mathbf{B})\mathbf{A}^{-it}\otimes \log(\rho)\bigr) (\mathbf{A}^{it}\otimes \rho^{it})\\
&=
Z^* \exp\bigl( -i  \log(e^{t} \mathbf{B})\otimes \log(\rho)\bigr) (\mathbf{A}^{it}\otimes \rho^{it})=(\I\otimes \rho^{-it})(\mathbf{A}^{it}\otimes \rho^{it})=(\mathbf{A}\otimes \I)^{it}
\end{split}\]
for $t\in \RR$. Furthermore, using that $\rho$ is group-like we have
\[
(\id\otimes \Delta)(Z_{[23]}^*\vv_{[13]})=
Z_{[23]}^*Z_{[24]}^*
\vv_{[13]}\vv_{[14]}
=
Z_{[23]}^*\vv_{[13]}Z_{[24]}^*
\vv_{[14]},
\]
i.e.~$Z_{[23]}^*\vv_{[13]}$ is a right representation of $\GG$ on $\LL^2(\GG)\otimes\LL^2(\RR)$, consequently $Z^*_{[34]}\vv_{[24]}$ is a right representation of $\GG$ on $\msf{H}\otimes\LL^2(\GG)\otimes \LL^2(\RR)$. We claim that
\begin{equation}\label{eq15}
q_{Y\otimes \mathbf{A}\otimes\I}\bigl( \vv^*_{[24]}Z_{[34]}(\zeta\otimes \eta) \bigr)=
q_{Y\otimes \mathbf{A}}(\zeta) \|\eta\|^2
\end{equation}
for $\zeta\in \Dom(q_{Y\otimes \mathbf{A}}),\eta\in \LL^2(\GG)$. As a preliminary step towards \eqref{eq15}, we first claim
\begin{equation}\label{eq20}
q_{Y\otimes \mathbf{A}\otimes\I}\bigl(\vv^*_{[24]} (\xi\otimes \kappa)\bigr)=
q_Y(\xi) q_{\mathbf{A}\otimes \rho}(\kappa)
\end{equation}
for $\xi\in \Dom(q_Y), \kappa\in \Dom(\mathbf{A}^{1/2})\odot \Dom(\rho^{1/2})$. To see this, write $\kappa=\sum_{n=1}^{N}\kappa^1_n\otimes \kappa^2_n$ and calculate using \eqref{eq22} (extended to sesquilinear forms)
\[\begin{split}
&\quad\;
q_Y(\xi) q_{\mathbf{A}\otimes \rho}(\kappa)=
\sum_{n,n'=1}^{N}
q_Y(\xi) q_{\mathbf{A}}(\kappa^1_n,\kappa^1_{n'})q_{\rho}(\kappa^2_n,\kappa^2_{n'})\\
&=
\sum_{n,n'=1}^{N}
q_{\mathbf{A}}(\kappa^1_n,\kappa^1_{n'})q_{Y\otimes\I}\bigl(\vv^*_{[23]}(\xi\otimes\kappa^2_n),
\vv^*_{[23]}(\xi\otimes\kappa^2_{n'})\bigr)\\
&=
\sum_{n,n'=1}^{N}
q_{Y\otimes\mathbf{A}\otimes \I}\bigl(\vv^*_{[24]}(\xi\otimes\kappa^1_n\otimes \kappa^2_n),
\vv^*_{[24]}(\xi\otimes\kappa^1_{n'}\otimes \kappa^2_{n'})\bigr)=q_{Y\otimes\mathbf{A}\otimes\I}\bigl(\vv^*_{[24]}(\xi\otimes\kappa)\bigr),
\end{split}\]
which establishes \eqref{eq20}.

Next we prove that equality \eqref{eq20} holds for all $\xi\in \Dom(q_Y), \kappa\in \Dom(\mathbf{A}^{1/2}\otimes \rho^{1/2})$. By \cite[Example 1.13, page 315]{Kato}, forms $q_{Y\otimes\mathbf{A}\otimes\I}, q_{\mathbf{A}\otimes\rho}$ are closed\footnote{For a von Neumann algebra $\operatorname{Q}\subseteq \B(\msf{K})$ and $T=A_T+\infty (\I-p_T)\in \oon{Q}^{+,\ext}$, sesquilinear form $q_T$ equals the form from \cite{Kato} for the closed operator $A_T^{1/2}$ with $\Dom(q_T)=\Dom(A_T^{1/2})\subseteq \oon{Im}(p_T)\subseteq\msf{K}$.}. Fix $\xi\in \Dom(Y^{1/2}),\kappa\in \Dom(\mathbf{A}^{1/2}\otimes \rho^{1/2})$. Take $(\kappa_n)_{n\in\NN}$ in $\Dom(\mathbf{A}^{1/2})\odot \Dom(\rho^{1/2})$ such that $\kappa_n\xrightarrow[n\to\infty]{q_{\mathbf{A}\otimes \rho}} \kappa$, which means $\kappa_n\xrightarrow[n\to\infty]{}\kappa$ and $q_{\mathbf{A}\otimes \rho}(\kappa_n-\kappa_m)\xrightarrow[n,m\to\infty]{}0$ \cite[Page 313]{Kato}. Such sequence exists, as $\Dom(\mathbf{A}^{1/2})\odot \Dom(\rho^{1/2})$ is a core for $\mathbf{A}^{1/2}\otimes \rho^{1/2}$. Since $\vv^*_{[24]}(\xi\otimes\kappa_n)\in \Dom(q_{Y\otimes\mathbf{A}\otimes\I})$ by \eqref{eq20}, $\vv^*_{[24]}(\xi\otimes \kappa_n)\xrightarrow[n\to\infty]{}\vv^*_{[24]}(\xi\otimes \kappa)$ and
\begin{equation}\label{eq24}
q_{Y\otimes \mathbf{A}\otimes \I}(
\vv^*_{[24]}(\xi\otimes \kappa_n) - 
\vv^*_{[24]}(\xi\otimes \kappa_m)
)=
q_Y(\xi)q_{\mathbf{A}\otimes\rho}(\kappa_n-\kappa_m)\xrightarrow[n,m\to\infty]{}0,
\end{equation}
we have $\vv^*_{[24]}(\xi\otimes \kappa)\in \Dom(q_{Y\otimes \mathbf{A}\otimes \I})$ and $\vv^*_{[24]}(\xi\otimes \kappa_n)\xrightarrow[n\to\infty]{q_{Y\otimes\mathbf{A}\otimes\I}} \vv^*_{[24]}(\xi\otimes \kappa)$. Using \cite[Theorem 1.12, Page 314]{Kato} we conclude that 
\[
q_Y(\xi)q_{\mathbf{A}\otimes \rho}(\kappa)=
\lim_{n\to\infty}
q_Y(\xi)q_{\mathbf{A}\otimes \rho}(\kappa_n)=
\lim_{n\to\infty} q_{Y\otimes \mathbf{A}\otimes \I}(\vv^*_{[24]} (\xi\otimes \kappa_n))=
q_{Y\otimes \mathbf{A}\otimes \I}(\vv^*_{[24]} (\xi\otimes \kappa)),
\]
which establishes \eqref{eq20} for $\kappa$ in the larger domain $\Dom(q_{\mathbf{A}\otimes\rho})$. Immediately from \eqref{eq14} follows $\mathbf{A}^{1/2}\otimes\rho^{1/2}=Z(\mathbf{A}^{1/2}\otimes\I)Z^*$ and consequently
\begin{equation}\label{eq23}
q_{Y\otimes\mathbf{A}\otimes\I}\bigl(\vv^*_{[24]}Z_{[34]}(\xi\otimes \varkappa)\bigr)=
q_Y(\xi)q_{\mathbf{A}\otimes\I}(\varkappa)
\end{equation}
for $\xi\in \Dom(q_Y), \varkappa\in \Dom(q_{\mathbf{A}\otimes\I})$. Take now $\zeta=\sum_{m=1}^{M}\zeta^1_m\otimes\zeta^2_m\in \Dom(q_Y)\odot \Dom(q_\mathbf{A}), \eta\in \LL^2(\GG)$. Using \eqref{eq23} we obtain \eqref{eq15} for this choice of vectors:
\[\begin{split}
&\quad\;
q_{Y\otimes\mathbf{A}\otimes\I}\bigl(\vv^*_{[24]}Z_{[34]}(\zeta\otimes\eta)\bigr)\\
&=
\sum_{m,m'=1}^{M} 
q_{Y\otimes\mathbf{A}\otimes\I}\bigl(\vv^*_{[24]}Z_{[34]}(\zeta^1_m\otimes\zeta^2_m\otimes\eta),\!
\vv^*_{[24]}Z_{[34]}(\zeta^1_{m'}\otimes\zeta^2_{m'}\otimes\eta)
\bigr)\\
&=
\sum_{m,m'=1}^{M}
q_Y(\zeta^1_m,\zeta^1_{m'})q_{\mathbf{A}\otimes\I}(\zeta^2_m\otimes \eta, 
\zeta^2_{m'}\otimes\eta)=
q_{Y\otimes\mathbf{A}}(\zeta)\|\eta\|^2.
\end{split}\]
Finally, take arbitrary $\zeta\in \Dom(q_{Y\otimes\mathbf{A}}),\eta\in \LL^2(\GG)$ and choose $\zeta_n\in \Dom(q_Y)\odot \Dom(q_{\mathbf{A}})$ with $\zeta_n\xrightarrow[n\to\infty]{q_{Y\otimes\mathbf{A}}}\zeta$. Arguing as in \eqref{eq24} we see $\vv^*_{[24]}Z_{[34]} (\zeta_n\otimes\eta)\xrightarrow[n\to\infty]{q_{Y\otimes\mathbf{A}\otimes\I}}\vv^*_{[24]}Z_{[34]} (\zeta\otimes\eta)$ and the general case of \eqref{eq15} holds by approximation.

Now we are in the situation of Lemma \ref{lemma1} (with von Neumann algebra $\OnP\bar\otimes\LL^{\infty}(\GG)\bar\otimes \B(\LL^2(\RR))$ and right representation $Z^*_{[34]}\vv_{[24]}$ of $\GG$ on $\msf{H}\otimes\LL^2(\GG)\otimes\LL^2(\RR)$). Using it and \eqref{eq14} we obtain
\[
Y\otimes\mathbf{A}\otimes \I
=
\vv_{[24]}^*Z_{[34]}(Y\otimes \mathbf{A}\otimes\I)
Z^*_{[34]}\vv_{[24]}=
\vv_{[24]}^* (Y\otimes \mathbf{A}\otimes\rho)
\vv_{[24]}.
\]
Clearly this implies
\begin{equation}\label{eq25}
Y\otimes\I=\vv^*_{[23]}(Y\otimes\rho)\vv_{[23]}.
\end{equation}
Uniqueness of spectral decomposition \eqref{eq42} gives 
\begin{equation}\label{eq26}
p_Y\otimes\I=\vv^*_{[23]}(p_Y\otimes\I)\vv_{[23]},\quad
A_Y\otimes\I=
\vv^*_{[23]}(A_Y\otimes\rho)\vv_{[23]}\;\textnormal{ on }\;\oon{Im}(p_Y)\otimes \LL^2(\GG).
\end{equation}
The second equation implies $1_{\{0\}}(A_Y)\otimes \I=\vv^*_{[23]}(1_{\{0\}}(A_Y)\otimes \I)\vv_{[23]}$.

Since $Y\in( \OnP\bar\otimes\LL^{\infty}(\GG))^{+,\ext}$ and $\LL^{\infty}(\GG)\cap\LL^{\infty}(\whG)=\CC\I$ \cite[Page 22]{VanDaeleVN}, we obtain $p_Y=e_{\oon{fin}}\otimes \I$ and $1_{\{0\}}(A_Y)=e_0\otimes \I$ for projections $e_{fin}\in \OnP, e_0\in e_{\oon{fin}}\OnP e_{\oon{fin}}$.

Consider now $B_Y=A_Y\rest_{\oon{Im}(e_{\oon{fin}}-e_0)\otimes\LL^2(\GG)}$, which is a strictly positive, self-adjoint operator. By \eqref{eq26} we still have $B_Y\otimes \I=\vv^*_{[23]} (B_Y\otimes\rho)\vv_{[23]}$ on $\oon{Im}(e_{\oon{fin}}-e_0)\otimes \LL^2(\GG)^{\otimes 2}$. Since $B_Y$ is strictly positive, self-adjoint we can apply imaginary powers:
\[
\vv_{[23]}(B_Y^{it}(\I\otimes\rho^{-it})\otimes\I)\vv_{[23]}^*=
(B_Y^{it}\otimes\rho^{it})(\I\otimes\rho^{-it}\otimes \rho^{-it})=
B^{it}(\I\otimes\rho^{-it})\otimes \I
\]
for $t\in\RR$, which implies $B_Y^{it}(\I\otimes\rho^{-it})\in (e_{\oon{fin}}-e_0)\OnP(e_{\oon{fin}}-e_0)\bar\otimes \CC\I$. Standard reasoning and application of Stone's theorem give strictly positive, self-adjoint operator $C_{2}$ affiliated with $(e_{\oon{fin}}-e_0)\OnP(e_{\oon{fin}}-e_0)$ such that $B_Y=C_2\otimes \rho$. Extend by $0$ operator $C_2$ to self-adjoint, positive operator $C_{\oon{fin}}$ affiliated with $e_{\oon{fin}}\OnP e_{\oon{fin}}$; then $A_Y=C_{\oon{fin}}\otimes\rho$. Finally, define $C\in \OnP^{+,\ext}$ by $C=C_{\oon{fin}}+\infty(\I-e_{\oon{fin}})$ so that $Y=C\otimes \rho$ as claimed. Uniqueness of $C$ is clear.
\end{proof}

With these preliminaries, we can prove Theorem \ref{thm2}. Recall that in Lemma \ref{lemma4} we have introduced a one parameter group of automorphisms $\kappa_t=\tau_t \Ad(\rho^{-it})\in \Aut(\LL^{\infty}(\GG))\;(t\in\RR)$.

\begin{proof}[Proof of Theorem \ref{thm2}]
It is clear that $\eqref{RI2}\Rightarrow \eqref{RI1}$ and $\eqref{RI4}\Rightarrow \eqref{RI3}$. To see $\eqref{RI3}\Rightarrow\eqref{RI2}$ take $\OnP=\CC$ and apply vector functional $\omega_\xi$ with $\xi\in \Dom(\rho^{1/2})$ to both sides of \eqref{RI3}. Consequently it is enough to show $\eqref{RI1}\Rightarrow\eqref{RI2}\Rightarrow \eqref{RI3}\Rightarrow \eqref{RI4}$.

Assume condition \eqref{RI1}. Fix von Neumann algebra $\OnP\subseteq \B(\LL^2(\OnP))$, $X\in (\OnP\bar\otimes \N)^{+}$ and consider element $(\id\otimes\id\otimes\theta)(\id\otimes\alpha)(X)\in (\OnP\bar\otimes \LL^{\infty}(\GG))^{+,\ext}$. For a vector $\xi\in \Dom(q_{(\id\otimes\id\otimes\theta)(\id\otimes\alpha)(X)})$ we have
\[
\theta\bigl(
(\omega_\xi\otimes\id)(\id\otimes \alpha)(X)\bigr)=
\omega_\xi\bigl((\id\otimes\id\otimes\theta)(\id\otimes \alpha)(X)\bigr)=
q_{(\id\otimes\id\otimes\theta)(\id\otimes\alpha)(X)}(\xi)<+\infty,
\]
hence $(\omega_\xi\otimes\id)(\id\otimes\alpha)(X)\in \mf{M}^+_{\theta}$ and for $\eta\in \Dom(\rho^{1/2})$
\[\begin{split}
&\quad\;
q_{(\id\otimes\Delta)  ((\id\otimes \id\otimes\theta)(\id\otimes \alpha)(X) )}(\xi\otimes \eta)=
(\omega_{\xi}\otimes\omega_{\eta})
\bigl(
(\id\otimes\Delta) 
( (\id\otimes \id\otimes\theta)(\id\otimes \alpha)(X) )
\bigr)\\
&=
\omega_{\eta}\bigl(
(\id\otimes\theta)
\alpha\bigl(
(\omega_\xi\otimes\id)
(\id\otimes\alpha)(X) 
\bigr)\bigr)=
\|\rho^{1/2}\eta\|^2 
\theta\bigl((\omega_\xi\otimes\id)(\id\otimes\alpha)(X)\bigr)\\
&=\|\rho^{1/2}\eta\|^2 
\omega_\xi\bigl((\id\otimes\id\otimes\theta )(\id\otimes\alpha)(X)\bigr)=
\|\rho^{1/2}\eta\|^2 q_{(\id\otimes\id\otimes\theta)(\id\otimes\alpha)(X)}(\xi),
\end{split}
\]
using the assumed relative $\rho$-invariance of $\theta$ on integrable elements \eqref{RI1}. Lemma \ref{lemma2} applies and we obtain
\begin{equation}\label{eq4}
(\id\otimes \id\otimes\theta)(\id\otimes \alpha)(X)=C_X\otimes \rho\quad\textnormal{in}\quad(\OnP\bar\otimes \LL^{\infty}(\GG))^{+,\ext}
\end{equation}
for a unique $C_X\in \OnP^{+,\ext}$.\\

First we establish $\eqref{RI1}\Rightarrow \eqref{RI2}$; specify $\OnP=\CC$, then $C_X\in \left[0,\infty\right]$. Take a vector $\zeta\in \Dom(\rho^{1/2})$ with $\|\rho^{1/2}\zeta\|=1$. Applying $\omega_\zeta$ to both sides of \eqref{eq4} shows that
\begin{equation}\label{eq5}
\Omega\colon \N^+\ni X\mapsto C_X=
 \theta\bigl(
( \omega_\zeta\otimes\id)\alpha(X)\bigr)\in \left[0,\infty\right]
 \end{equation}
 is a normal weight on $\N$ which does not depend on the choice of $\zeta$. For $X\in\mf{M}_{\theta}^{+}$ we have
\begin{equation}\label{eq6}
\Omega(X)=
\theta\bigl((\omega_\zeta\otimes\id)\alpha(X)\bigr)
=\|\rho^{1/2}\zeta\|^2
\theta(X)=\theta(X)<+\infty
\end{equation}
using \eqref{RI1}. In particular 
$X\in \mf{M}^+_{\Omega}$ and $\Omega$ is semifinite. Finally, $\Omega$ is faithful. Indeed, this follows from \eqref{eq4} as both $\alpha$ and $ \id\otimes \theta$ are faithful. We need to show that $\Omega = \theta$ as n.s.f.~weights on $\N$.

Take $t\in\RR,X\in\N^+$ and calculate using Lemma \ref{lemma4}
\[\begin{split}
\Omega(\sigma^{\theta}_t(X))&=
\theta\bigl((\omega_\zeta\otimes \id)\bigl(\alpha( \sigma^{\theta}_t(X))\bigr)\bigr)=
\theta\bigl(
(\omega_\zeta\otimes \id)
(\kappa_t\otimes \sigma^\theta_t)
\alpha(X)\bigr)\\
&=
\theta\bigl(
(\omega_{\rho^{it} P^{-it}\zeta}  \otimes \id)
\alpha(X)\bigr)
=
C_X\, \omega_{\rho^{it}P^{-it}\zeta} (\rho)=
C_X=\Omega(X).
\end{split}\]
As $X$ is arbitrary, it follows that $\Omega \sigma^{ \theta}_t=\Omega$, i.e.~weights $\Omega, \theta$ commute \cite[Theorem 4.10]{Stratila}. Since these weights are equal on $\mf{M}_{\theta}^+$ \eqref{eq6}, \cite[Theorem 6.2]{Stratila} applies and we obtain $\Omega=\theta$. This means $C_X=\Omega(X)=\theta(X)$ for $X\in \N^+$ and equality \eqref{eq4} gives $(\id\otimes \theta)\alpha(X)=\theta(X)\rho$. Applying $\omega_\xi$ with $\xi\in \Dom(\rho^{1/2})$ ends the proof of \eqref{RI2}.\\

$\eqref{RI2}\Rightarrow\eqref{RI3}$; consider again general $\OnP, X\in (\OnP\bar\otimes\N)^+$. Take $\zeta\in \Dom(\rho^{1/2}), \|\rho^{1/2}\zeta\|=1, \omega\in \OnP^+_*$ and apply $\omega\otimes \omega_\zeta$ to \eqref{eq4}:
\[
\omega(C_X)=(\omega\otimes \omega_\zeta)(C_X\otimes \rho)=
\theta \bigl((\omega_\zeta\otimes\id)\alpha\bigl((\omega\otimes\id)X\bigr)\bigr)=
\|\rho^{1/2}\zeta\|^2 \theta\bigl((\omega\otimes \id)X\bigr)=
\omega\bigl((\id\otimes \theta)X\bigr).
\]
Note that we use \eqref{RI2} in third equality, as we do not know whether $(\omega\otimes\id)(X)$ is in $\mf{M}_\theta^+$. The above equation holds for all $\omega\in\OnP^+_*$, hence we conclude $C_X=(\id\otimes \theta)X$. Inserting this equality into \eqref{eq4} gives \eqref{RI3}.\\

$\eqref{RI3}\Rightarrow\eqref{RI4}$; Take a von Neumann algebra $\OnP$ and $X\in (\OnP\bar\otimes \N)^{+,\ext}$. We can find a sequence $X_n\in (\OnP\bar\otimes \N)^+\,(n\in\NN)$ such that $X$ is the increasing limit of $(X_n)_{n\in\NN}$ \cite[Section 11.4]{Stratila}. Then using \eqref{RI3} and \cite[Proposition 11.4]{Stratila} we get
\[\begin{split}
&\quad\;
\omega\bigl( (\id\otimes\id\otimes\theta)(\id\otimes\alpha)(X)\bigr)=
\sup_{n\in\NN}
\omega\bigl( (\id\otimes\id\otimes\theta)(\id\otimes\alpha)(X_n)\bigr)\\
&=
\sup_{n\in\NN}\omega\bigl( (\id\otimes\theta)(X_n)\otimes\rho\bigr)=
\omega\bigl( (\id\otimes\theta)(X)\otimes\rho\bigr)
\end{split}\]
for every $\omega\in (\OnP\bar\otimes \LL^{\infty}(\GG))_*^+$, which is exactly \eqref{RI4}.
\end{proof}

\section{Universal positive cone condition}\label{sec:univ}

Quantum groups are studied at two levels: reduced and universal. Since there is the canonical reducing map $\pi_{\GG}\colon \mrm{C}_0^u(\GG)\rightarrow \mrm{C}_0(\GG)$, it is often easy to transport structure or properties from universal to reduced level. The reverse direction is often less obvious. Nonetheless, some objects admit canonical lifts: for example one can lift representations and bicharacters \cite{Kustermans} or von Neumann algebraic actions \cite[Section 4.3]{DeCommerKrajczok}. Being able to pass between these two pictures is often technically convenient. We will use such passage to prove Corollary \ref{cor1}, which intuitively states that standard representation composed with quantum group homomorphism is standard. This result is motivated by a joint work of Kenny De Commer and the author on braided locally compact quantum groups \cite{DeCommerKrajczokBLCQG}.

Our first task is to formulate and establish a version of conditions from Theorem \ref{thm1} at the universal level. Implementation of action and the condition involving modular conjugations will not cause much problems, and we will prove it in Proposition \ref{prop3}. The universal version of condition $(3)$, concerning positive cone, requires more delicate handling. We state the relevant definition in Definition \ref{def2} and prove it in Theorem \ref{thm5}. 

Recall the notions of universal lift of representation $\mathds{U}$ and action $\alpha^{u}$ from the Preliminaries. We begin with a formula for the universal lift in the relatively invariant case, similar to the expression at the reduced level (Example \ref{ex1}). This result is likely to be of independent interest.

\begin{lemma} \label{lemma5}
Let $\alpha\colon \GG\curvearrowright \N$ be an action, $\rho\,\eta\LL^{\infty}(\GG)$ a strictly positive, self-adjoint group-like operator and $\theta\in \oon{W}_{\nsf}(\N)$ a relatively $\rho$-invariant weight. Let $\mrm{U}^\alpha$ be the corresponding standard implementation of $\alpha$ with the universal lift $\mathds{U}^\alpha$. 
\begin{itemize}
\item The space $\mf{N}_{\theta,0}$ spanned by elements of the form $ (\omega_{\xi,\eta}\otimes\id)\alpha (n)$ with $n\in\mf{N}_\theta$ analytic for $\sigma^\theta$ and $\xi,\eta\in\bigcap_{k\in\ZZ} \Dom(P^k)\cap\Dom(\rho^k)$, forms a $\ssots\times \|\cdot\|$ core for $\Lambda_\theta$.

\item Consider functional of the form $\omega_{\xi,\zeta}\pi\in \mrm{C}_0^u(\GG)^*$, where $\pi\colon \mrm{C}_0^u(\GG)\rightarrow \B(\msf{K})$ is a non-degenerate representation and $\xi\in \msf{K}, \zeta\in \Dom (\pi(\rho^u) ^{-1/2})$. For $n\in\mf{N}_\theta$ we have 
\[
(\omega_{\xi,\pi(\rho^u)^{-1/2}\zeta }\pi\otimes \id)\alpha^u(n)\in \mf{N}_\theta
\]
and
\begin{equation}\label{eq33}
(\omega_{\xi,\zeta}\pi\otimes \id)(\mathds{U}^{\alpha *})\Lambda_\theta(n)=\Lambda_\theta ((\omega_{\xi,\pi(\rho^u)^{-1/2}\zeta }\pi\otimes \id)\alpha^u(n)).
\end{equation}
\end{itemize}
\end{lemma}

\begin{proof}[Sketch of proof]
Consider $\Lambda_\theta\rest_{\mf{N}_{\theta,0}}$. This map is $\ssots\times\|\cdot\|$ closable, let $\mc{D}$ be the domain of its closure. We claim that $\mc{D}$ has the following properties:
\begin{enumerate}
\item $\sigma^\theta_t(\mc{D})=\mc{D}$ for $t\in\RR$,
\item $\mc{D}$ is a $\ssots$-dense left ideal in $\N$.
\end{enumerate}
Lemma \ref{lemma4} gives $\sigma^\theta_t(\mf{N}_{\theta,0})=\mf{N}_{\theta,0}$ and $\sigma^\theta_t(\mc{D})=\mc{D}$ follows by approximation. Take $n\in\N$. Since Podleś condition is automatic in the von Neumann algebraic setting \cite[Proposition 2.6]{QGProj1} (see also \cite[Section 3.2]{CrannDeRoKrajczok}), we can find a net $(n_i)_{i\in I}$ in $\lin\{(\omega\otimes\id)\alpha(m)\mid \omega\in \LL^1(\GG),m\in\N\}$ such that $n_i\xrightarrow[i\in I]{\ssots} n$. By a standard approximation argument, we can in fact assume that $n_i\in \mf{N}_{\theta,0}$; this shows $\ssots$-density of $\mf{N}_{\theta,0}$ in $\N$ and in turn the density of $\mc{D}$. Finally, we need to argue that $\mc{D}$ is a left ideal. Take $x\in \N, n\in\mc{D}$ and choose $(n_i)_{i\in  I}$ in $\mf{N}_{\theta,0}$ which approximates $n$ in the graph topology. It is enough to show $x n_i\in \mc{D}$, as then $xn\in \mc{D}$ follows by the fact that $\Lambda_{\theta}\rest_{\mc{D}}$ is closed. Write $n_i=\sum_{k=1}^{K} (\omega_{\xi_k,\eta_k}\otimes \id)\alpha(m_k)\in \mf{N}_{\theta,0}$ and choose $x_{l,j}\in \N,y_{l,j}\in \LL^{\infty}(\GG)$ such that
\[
\sum_{l=1}^{L_j} (y_{l,j}\otimes\I)\alpha(x_{l,j})\xrightarrow[j\in J]{\ssots}\I\otimes x
\]
(using again Podleś condition). By further approximation, we can assume that $y_{l,j}$'s are analytic for $\tau,(\Ad(\rho^{it}))_{t\in\RR}$ and $x_{l,j}$'s are analytic for $\sigma^\theta$. We have
\[
\mf{N}_{\theta,0}\ni \sum_{k=1}^K \sum_{l=1}^{L_j}
(\omega_{y_{l,j}^* \xi_k,\eta_k}\otimes\id)\alpha(x_{l,j} m_k)\xrightarrow[j\in J]{\ssots}
\sum_{k=1}^{K}(\omega_{\xi_k,\eta_k}\otimes\id)\bigl((\I\otimes x)\alpha(m_k)\bigr)=x n_i
\]
and using Example \ref{ex1} 
\[\begin{split}
&\quad\;
\Lambda_\theta\bigl(
\sum_{k=1}^K \sum_{l=1}^{L_j}
(\omega_{y_{l,j}^* \xi_k,\eta_k}\otimes\id)\alpha(x_{l,j} m_k)
\bigr)=
\sum_{k=1}^K \sum_{l=1}^{L_j}
(\omega_{y_{l,j}^* \xi_k,\rho^{1/2}\eta_k}\otimes\id)(\mrm{U}^{\alpha *})
\Lambda_\theta(
 x_{l,j} m_k)\\
 &=
\sum_{k=1}^K 
(\omega_{  \xi_k,\rho^{1/2}\eta_k}\otimes\id)
\bigl(
\sum_{l=1}^{L_j}
(y_{l,j}\otimes\I)\alpha(x_{l,j})\mrm{U}^{\alpha *}\bigr)
\Lambda_\theta(m_k)\xrightarrow[j\in J]{}
x
\sum_{k=1}^K 
(\omega_{  \xi_k,\rho^{1/2}\eta_k}\otimes\id)
(\mrm{U}^{\alpha *})
\Lambda_\theta(m_k),
\end{split}\]
which shows that $x n_i\in \mc{D}$ and $\mc{D}$ is a left ideal. Now, take any $n\in \mf{N}_\theta$, choose a bounded net $(e_\lambda)_{\lambda\in \Lambda}$ in $\mc{D}$ with $e_\lambda\xrightarrow[\lambda\in \Lambda]{\ssots}\I$ (Kaplansky's density) and set $e_{\lambda,0}=\tfrac{1}{\sqrt{\pi}}\int_{\RR}e^{-s^2} \sigma^\theta_s(e_\lambda)\md s$. By Hille's theorem $e_{\lambda,0}\in \mc{D}$ (see e.g.~\cite[Lemma 14.10]{KrajczokCoamenability}). We have $n e_{\lambda,0}\in \mc{D}$ and $\Lambda_\theta(n e_{\lambda,0})=J_\theta \sigma^\theta_{i/2}(e_{\lambda,0})^* J_\theta \Lambda_\theta(n)\xrightarrow[\lambda\in \Lambda]{}\Lambda_\theta(n)$, which shows $n\in\mc{D}$ and ends the proof of the first bullet point.

To prove the second bullet point, we need the following extension of the formula for $\mrm{U}^{\alpha}$ (Example \ref{ex1}): take $n\in\mf{N}_\theta$ and $\omega\in \LL^1(\GG)$ such that $\rho^{-1/2}\omega\in \LL^1(\GG)$. By this statement we mean that there is $\nu\in \LL^1(\GG)$ such that $\omega(x\rho^{-1/2})=\nu(x)$ for $x\in \LL^{\infty}(\GG)$ with $x\rho^{-1/2}$ bounded on $\Dom(\rho^{-1/2})$. Then
\begin{equation}\label{eq32}
(\rho^{-1/2}\omega\otimes\id)\alpha(n)\in\mf{N}_\theta\quad\textnormal{and}\quad
(\omega\otimes\id)(\mrm{U}^{\alpha *})\Lambda_\theta(n)=\Lambda_\theta((\rho^{-1/2}\omega\otimes\id)\alpha(n)).
\end{equation}
One can prove \eqref{eq32} by approximation, using that $\Lambda_\theta$ is closed.

Now, take $\nu=\omega_{\xi,\eta}\pi\in \mrm{C}_0^u(\GG)^{*}$ as in the claim and $n=(\omega_{\zeta,\rho^{-1/2}\chi}\otimes\id)\alpha(m)\in \mf{N}_{\theta,0}$. Then
\[\begin{split}
&\quad\;
(\nu\otimes\id)(\mathds{U}^{\alpha *})\Lambda_\theta(n)=
(\nu\otimes\id)(\mathds{U}^{\alpha *})
(\omega_{\zeta,\chi}\otimes\id)(\mrm{U}^{\alpha *})
\Lambda_\theta(m)\\
&=
(\omega_{\zeta,\chi}\star\nu\otimes\id)(\mrm{U}^{\alpha *})\Lambda_\theta(m)=
\Lambda_\theta\bigl( (\rho^{-1/2} (\omega_{\zeta,\chi}\star\nu)\otimes\id)\alpha(m)\bigr)
\end{split}\]
using \eqref{eq32}, with $ \omega_{\zeta,\chi}\star\nu=(\omega_{\zeta,\chi}\otimes\nu)\Delta^{r,u,\vN}\in \LL^1(\GG)$ (see \cite[Section 4.3]{DeCommerKrajczok}). Since $\rho^{-1/2} (\omega_{\zeta,\chi}\star\nu)=\omega_{\zeta,\rho^{-1/2}\chi}\star \omega_{\xi,\pi(\rho^u)^{-1/2}\eta}\pi \in \LL^1(\GG)$, this further equals
\[
\Lambda_\theta\bigl(
(\omega_{\xi,\pi(\rho^u)^{-1/2}\eta}\pi\otimes\id)\alpha^u\bigl(
 ( \omega_{\zeta,\rho^{-1/2}\chi}\otimes\id)\alpha(m)
\bigr)\bigr)=
\Lambda_\theta\bigl(
(\omega_{\xi,\pi(\rho^u)^{-1/2}\eta}\pi\otimes\id)\alpha^u(
n)\bigr)
\]
and proves \eqref{eq33} for $n\in\mf{N}_{\theta,0}$. General case follows by approximation, using that $\mf{N}_{\theta,0}$ is a $\ssots\times\|\cdot\|$ core for $\Lambda_\theta$.
\end{proof}

As a corollary, we deduce a universal version of the relative invariance condition \eqref{RI1}.

\begin{proposition}\label{prop2}
Let $\alpha\colon \GG\curvearrowright \N$ be an action, $\rho\,\eta\LL^{\infty}(\GG)$ a strictly positive, self-adjoint group-like operator and $\theta\in \oon{W}_{\nsf}(\N)$ a relatively $\rho$-invariant weight. We have
\[
\theta\bigl( (\omega_\xi\pi\otimes\id)\alpha^{u}(n)\bigr)=\|\pi(\rho^u)^{1/2}\xi\|^2 \theta(n)
\]
for $n\in \mf{M}_{\theta}^+$, non-degenerate representation $\pi\colon \mrm{C}_0^u(\GG)\rightarrow \B(\msf{K})$ and $\xi\in \Dom(\pi(\rho^u)^{1/2} )\subseteq \msf{K}$.
\end{proposition}

\begin{proof}
Denote by $\mrm{U}^\alpha$ the standard implementation of $\alpha$ with respect to weight $\theta$, with universal lift $\mathds{U}^\alpha$. Let $\{\xi_i\}_{i\in I}$ be an orthonormal basis in $\msf{K}$. We calculate using Lemma \ref{lemma5} and normality of $\theta$:
\[\begin{split}
&\quad\;
\theta\bigl( (\omega_\xi\pi\otimes\id)\alpha^{u}(n)\bigr)=
\theta\bigl(\sum_{i\in  I}\bigl| (\omega_{\xi_i,\xi}\pi\otimes\id)\alpha^{u}(n^{1/2})\bigr|^2 \bigr)=
\sum_{i\in  I} \theta\bigl(\bigl| (\omega_{\xi_i,\xi}\pi\otimes\id)\alpha^{u}(n^{1/2})\bigr|^2 \bigr)\\
&=
\sum_{i\in I}\bigl\|
\Lambda_\theta\bigl(
(\omega_{\xi_i,\xi}\pi\otimes\id)\alpha^{u}(n^{1/2})
\bigr)\bigr\|^2=
\sum_{i\in I}\bigl\|
(\omega_{\xi_i,\pi(\rho^u )^{1/2}\xi}\pi \otimes\id)(\mathds{U}^{\alpha *}) \Lambda_\theta (n^{1/2})
\bigr\|^2 \\
&=
\bigl\la 
\Lambda_\theta(n^{1/2})\big|
\bigl(\,
\sum_{i\in I}
(\omega_{\pi(\rho^u )^{1/2}\xi,\xi_i}\pi \otimes\id)(\mathds{U}^{\alpha }) 
(\omega_{\xi_i,\pi(\rho^u )^{1/2}\xi}\pi \otimes\id)(\mathds{U}^{\alpha *}) 
\bigr)\Lambda_\theta(n^{1/2})\bigr\ra \!=\!
\|\pi(\rho^u)^{1/2}\xi\|^2 \, \theta(n).
\end{split}\]
\end{proof}

As explained at the beginning of this section, we want to express and prove conditions similar to these from Theorem \ref{thm1} at the universal level. Conditions $(1), (2)$ are not problematic.

\begin{proposition}\label{prop3}
Let $\N$ be a von Neumann algebra with a weight $\theta\in \oon{W}_{\nsf}(\N)$ and action $\alpha\colon \GG\curvearrowright \N$ of a locally compact quantum group. Let $\mrm{U}^{\alpha}$ be the corresponding standard implementation with lift $\mathds{U}^{\alpha}$. We have
\begin{enumerate}
\item $\mathds{U}^{\alpha *}(\I\otimes n)\mathds{U}^\alpha=\alpha^u(n)$ for $n\in\N$,
\item $(R^{u}\otimes j_\theta)(\mathds{U}^{\alpha})=\mathds{U}^\alpha$, where $j_\theta=J_\theta (\cdot)^* J_\theta$.
\end{enumerate}
\end{proposition}

\begin{proof}
Equality $\mathds{U}^{\alpha *}(\I\otimes n)\mathds{U}^\alpha=\alpha^u(n)$ is a definition of $\alpha^u$, see \cite[Section 4.3]{DeCommerKrajczok}. For the second point, write $\mrm{U}^{\alpha}=(\id\otimes \phi_\alpha)(\wW^{\GG})$, then $\mathds{U}^\alpha=(\id\otimes \phi_\alpha)(\WW^{\GG})$. From Theorem \ref{thm1} we know $j_\theta \phi_\alpha \wh{R}^u=\phi_\alpha$, hence equality $(R^{u}\otimes j_\theta)(\mathds{U}^{\alpha})=\mathds{U}^\alpha$ follows from $(R^u\otimes \wh{R}^u)(\WW^{\GG})=\WW^{\GG}$. 
\end{proof}

Condition $(3)$ in Theorem \ref{thm1} (involving positive cone) is more problematic, as one needs to put some analytical conditions on the relevant functionals, corresponding to the assumptions on vectors in Theorem \ref{thm1}. Such condition will be most convenient to state using representation and vector functionals.

\begin{definition}\label{def2}
Let $\N$ be a von Neumann algebra with weight $\theta\in \oon{W}_{\nsf}(\N)$ and $\mrm{U}\in \LL^{\infty}(\GG)\bar\otimes\B(\msf{H}_\theta)$ a unitary representation with the universal lift $\mathds{U}$. Let $n\in\NN, \omega\in \M_n(\mrm{C}_0^u(\GG)^*)$ be such that:
\begin{itemize}
\item there is a non-degenerate representation $\pi\colon \mrm{C}_0^u(\GG)\rightarrow \B(\msf{K})$ and vectors $\xi_1,\dotsc,\xi_n\in \msf{K}$ with $\omega=(\omega_{\xi_i,\xi_j}\pi)_{i,j=1}^{n}$,
\item there is a strictly positive, self-adjoint operator $A$ on $\msf{K}$ such that $\pi \tau^u_t=\Ad(A^{it})\pi$ for $t\in \RR$, 
\item it holds that $\xi_1,\dotsc,\xi_n\in \Dom(A^{1/4})\cap\Dom(A^{-1/4})$.
\end{itemize}
We say that $\mrm{U}$ satisfies \emph{universal positive cone condition for} $\omega$ if the matrix
\[
\bigl((\omega_{\xi_i,\xi_j}\pi \tau^u_{i/4}\otimes \id)(\mathds{U}^*)\bigr)_{i,j=1}^{n}\in \M_n(\B(\msf{H}_\theta))
\]
preserves the positive cone in the sense of Definition \ref{def1}.
\end{definition}

\begin{remark}\noindent
\begin{itemize}
\item The above condition is a generalisation of property $(3)$ in Theorem \ref{thm1}: if we take $\msf{K}=\LL^2(\GG), \pi=\pi_{\GG}$ and $A=P$, then the positive cone condition reduces to the statement of $(3)$ in Theorem \ref{thm1}.
\item For each $1\le i,j\le n$, functional $\omega_{\xi_i,\xi_j}\pi\tau^u_{i/4}$ is well defined in $\mrm{C}_0^u(\GG)^*$ as $\omega_{A^{-1/4}\xi_i,A^{1/4}\xi_j}\pi$.
\item The universal positive cone condition depends only on $\omega\in \M_n(\mrm{C}_0^u(\GG)^*)$, not on the choice of representation $\pi$ and vectors $\xi_1,\dotsc,\xi_n$.
\end{itemize}
\end{remark}

\begin{lemma}\label{lemma6}
Let $\alpha\colon\GG\curvearrowright \N$ be an action and $\theta_1,\theta_2\in \oon{W}_{nsf}(\N)$ weights with the corresponding standard implementations $\mrm{U}^{\alpha,1},\mrm{U}^{\alpha,2}$ of $\alpha$. Take $n\in\NN,\omega\in\M_n(\mrm{C}_0^u(\GG)^*)$ as in Definition \ref{def2}. If $\mrm{U}^{\alpha,1}$ satisfies the universal positive cone condition for $\omega$, then so does $\mrm{U}^{\alpha,2}$.
\end{lemma}

Proof of Lemma \ref{lemma6} is straightforward, hence will be skipped (cf.~Remark \ref{rem1}).\\

In Theorem \ref{thm5} we will establish the universal positive cone condition for all functionals satisfying the condition from  Definition \ref{def2}. As the proof uses the biduality theorem, we start with the relatively invariant case. 

\begin{proposition}\label{prop1}
Let $\N$ be a von Neumann algebra with weight $\theta\in \oon{W}_{\nsf}(\N)$, action $\alpha\colon \GG\curvearrowright \N$ and standard implementation $\mrm{U}^\alpha$. Assume that $\theta$ is relatively invariant with respect to a strictly positive, self-adjoint, group-like $\rho\,\eta\LL^{\infty}(\GG)$. Take $n\in\NN, \omega\in \M_n(\mrm{C}_0^u(\GG)^*)$ as in Definition \ref{def2}. Representation $\mrm{U}^{\alpha}$ satisfies the universal positive cone condition for $\omega$.
\end{proposition} 

\begin{proof}
Let $\xi_1,\dotsc\xi_n,\pi,A$ be objects from Definition \ref{def2}, in particular $\omega=(\omega_{\xi_i,\xi_j}\pi)_{i,j=1}^{n}$. Note first that $A$ and $\pi(\rho^u)$ strongly commute. Indeed,
\[
A^{it} \pi(\rho^u)^{is}=A^{it} \pi( \rho^{u is})=\pi (\tau^u_t(\rho^{u is})) A^{it}=\pi(\rho^u)^{is} A^{it}
\]
for all $s,t\in\RR$, using that $\rho^{u is}$ is a group-like unitary, hence it is invariant under scaling group (see e.g.~the discussion after Proposition 3.5 in \cite{KalantarNeufang}). We need to show 
\begin{equation}\label{eq34}
\bigl\la (\zeta'_{i,j})_{i,j=1}^{n}\big|
\bigl( (\omega_{\xi_i,\xi_j} \pi \tau^u_{i/4}\otimes\id)(\mathds{U}^{\alpha *}) \zeta_{i,j}\bigr)_{i,j=1}^{n}\bigr\ra \ge 0
\end{equation}
for $(\zeta'_{i,j})_{i,j=1}^{n},(\zeta_{i,j})_{i,j=1}^{n}\in \mf{P}_{\Tr_n\otimes \theta}$. By approximation, it is enough to assume that $\xi_i$'s live in a spectral subspace $1_{\left[\eps,\eps^{-1}\right]}(\pi(\rho^u)) 1_{[\eps,\eps^{-1}]}(A)\msf{K}$ for some $\eps>0$ and that $\zeta'_{i,j}=y_i J_\theta \Lambda_\theta(y_j),\zeta_{i,j}=x_i J_\theta \Lambda_\theta(x_j)$, where $x_i,y_i$'s and their adjoints are analytic for $\sigma^\theta$ and square integrable after applying $\sigma^\theta_z\,(z\in\CC)$ (Lemma \ref{lemma3}). With these reductions, we calculate
\[\begin{split}
&\quad\;
\bigl\la (\zeta'_{i,j})_{i,j=1}^{n}\big|
\bigl( (\omega_{\xi_i,\xi_j} \pi \tau^u_{i/4}\otimes\id)(\mathds{U}^{\alpha *}) \zeta_{i,j}\bigr)_{i,j=1}^{n}\bigr\ra\\
&=
\sum_{i,j=1}^{n}
\bigl\la y_i J_\theta \Lambda_\theta(y_j)| 
(\omega_{ A^{-1/4} \xi_i, A^{1/4}\xi_j}\pi\otimes\id)(\mathds{U}^{\alpha *}) x_i J_\theta \Lambda_\theta(x_j)\bigr\ra \\
&=
\sum_{i,j=1}^{n}
\bigl\la y_i J_\theta \Lambda_\theta(y_j)| 
(\omega_{ A^{-1/4} \xi_i, A^{1/4}\xi_j}\pi\otimes\id)(\alpha^u(x_i)\mathds{U}^{\alpha *}) J_\theta \Lambda_\theta(x_j)\bigr\ra \\
&=
\sum_{i,j=1}^{n}
\sum_{p\in \mc{P}}
\bigl\la \Lambda_\theta(y_i \sigma^\theta_{i/2}(y_j)^*)| 
(\omega_{ A^{-1/4} \xi_i, A^{1/4} \pi(\rho^u)^{-1/4} e_p}\pi\otimes\id)(\alpha^u(x_i))\\
&\qquad\qquad\qquad\qquad\qquad
\qquad\qquad
(\omega_{A^{-1/4} \pi(\rho^u)^{1/4}e_p,A^{1/4}\xi_j}\pi\otimes\id)(\mathds{U}^{\alpha *}) \Lambda_\theta(\sigma^\theta_{i/2}(x_j)^*)\bigr\ra 
\end{split}\]
where $\{e_p\}_{p\in \mc{P}}$ is an orthonormal basis in $\msf{K}$, chosen such that $e_p\in \bigcap_{k\in \ZZ} \Dom(A^k \pi(\rho^u)^{-k})$ and  $\|\sum_{p\in F} |A^{1/4} \pi(\rho^u)^{-1/4} e_p\ra \la 
A^{-1/4}  \pi(\rho^u)^{1/4}e_p|\|\le 2$ for all finite $F\subseteq \mc{P}$. One can obtain basis $\{e_p\}_{p\in \mc{P}}$ by choosing an orthonormal basis separately in each spectral subspace $1_{[2^k,2^{k+1})}(A^{1/4}\pi(\rho^u)^{-1/4})\msf{K}$ $(k\in\ZZ)$ and taking union, cf.~\cite[Lemma A.2]{MRWSemidirect}. Next, using lemmas \ref{lemma5}, \ref{lemma4} and \cite[Proposition 2.17]{Stratila}
\[\begin{split}
&\quad\;
\sum_{p\in \mc{P}}
\sum_{i,j=1}^{n}
\bigl\la \Lambda_\theta(y_i \sigma^\theta_{i/2}(y_j)^*)| 
(\omega_{ A^{-1/4} \xi_i, A^{1/4} \pi(\rho^u)^{-1/4} e_p}\pi\otimes\id)(\alpha^u(x_i))\\
&\qquad\qquad\qquad\qquad\qquad
\qquad\qquad
\Lambda_\theta\bigl(
(\omega_{A^{-1/4} \pi(\rho^u)^{1/4}e_p,\pi(\rho^u)^{-1/2}A^{1/4}\xi_j}\pi\otimes\id)\alpha^u\bigl(
 \sigma^\theta_{i/2}(x_j)^*\bigr)
\bigr)\bigr\ra \\
&=
\sum_{p\in \mc{P}}
\sum_{i,j=1}^{n}
\theta\bigl(
\sigma^\theta_{i/2}(y_j) y_i^* 
(\omega_{ A^{-1/4} \xi_i, A^{1/4} \pi(\rho^u)^{-1/4} e_p}\pi\otimes\id)(\alpha^u(x_i))\\
&\qquad\qquad\qquad\qquad\qquad
\qquad\qquad
(\omega_{A^{-1/4} \pi(\rho^u)^{1/4}e_p,\pi(\rho^u)^{-1/2}A^{1/4}\xi_j}\pi\otimes\id)\alpha^u\bigl(
 \sigma^\theta_{i/2}(x_j)^*\bigr)\;
\bigr)\\
&=
\sum_{p\in \mc{P}}
\sum_{i,j=1}^{n}
\theta\bigl(
\sigma^\theta_{i/2}(y_j) y_i^* 
(\omega_{ A^{-1/4} \xi_i, A^{1/4} \pi(\rho^u)^{-1/4} e_p}\pi\otimes\id)(\alpha^u(x_i))\\
&\qquad\qquad\qquad\qquad\qquad
\qquad\qquad
\sigma^\theta_{i/2}\bigl(
\,(\omega_{A^{-1/4}\xi_j,A^{1/4} \pi(\rho^u)^{-1/4}e_p}\pi\otimes\id)\alpha^u(x_j)\,\bigr)^*
\;\bigr)\\
&=
\sum_{p\in \mc{P}}
\theta\bigl(\bigl(\sum_{i=1}^{n}
 y_i ^* 
(\omega_{ A^{-1/4} \xi_i, A^{1/4} \pi(\rho^u)^{-1/4} e_p}\pi\otimes\id)(\alpha^u(x_i))\bigr)\\
&\qquad\qquad\qquad\qquad\qquad
\qquad\qquad
\sigma^\theta_{i/2}\bigl(\sum_{j=1}^{n} 
y_j^*
\,(\omega_{A^{-1/4}\xi_j,A^{1/4} \pi(\rho^u)^{-1/4}e_p}\pi\otimes\id)\alpha^u(x_j)\,\bigr)^*
\;\bigr)\ge 0.
\end{split}\]
The above inequality $\ge 0$ follows from the following general statement.
\[
\textnormal{If }z\in \mf{N}_\theta^*\cap \Dom(\sigma^\theta_{i/2})\textnormal{ and }
\sigma^\theta_{i/2}(z)\in \mf{N}_\theta^*, \textnormal{ then }\theta(z \sigma^\theta_{i/2}(z)^*)\ge 0,
\]
which for appropriate elements $z$ follows from $\theta(z \sigma^\theta_{i/2}(z)^*)=\la \Lambda_\theta(z^*) | \nabla_\theta^{1/2} \Lambda_\theta(z^*)\ra \ge 0$, and can be established in general by approximation.
\end{proof}

\begin{theorem}\label{thm5}
Let $\N$ be a von Neumann algebra with weight $\theta\in \oon{W}_{\nsf}(\N)$, action $\alpha\colon \GG\curvearrowright \N$ and the standard implementation $\mrm{U}^\alpha$ with respect to $\theta$. Take $n\in\NN, \omega\in \M_n(\mrm{C}_0^u(\GG)^*)$ as in Definition \ref{def2}. Representation $\mrm{U}^{\alpha}$ satisfies the universal positive cone condition for $\omega$.
\end{theorem}

\begin{proof}
Denote by $\xi_1,\dotsc,\xi_n, \pi, A$ the objects given by Definition \ref{def2}. The proof will be completed in two steps.

(Step 1) As recalled in Section \ref{sec:characterisation} and Example \ref{ex2.5}, one can define crossed product $\GG\ltimes\N\subseteq \B(\LL^2(\GG)\otimes \msf{H}_\theta)$ which carries a dual action $\alpha^{\ltimes}\colon \wh{\GG}^{\op}\curvearrowright \GG\ltimes \N$ and relatively $\delta_{\whG}$-invariant dual weight $\wh{\theta}\in \oon{W}_{\nsf}(\GG\ltimes\N)$. We can perform this operation again and obtain bidual action $\alpha^{\ltimes\ltimes}\colon {\GG^{\op}} '\curvearrowright \whG^{\op}\ltimes (\GG\ltimes \N)$ and bidual weight $\wh{\wh{\theta}}\in \oon{W}_{\nsf}(\whG^{\op}\ltimes(\GG\ltimes\N))$, which is relatively $\delta_{\GG'}$-invariant.

Consider now two actions $\beta,\gamma\colon \GG\curvearrowright \B(\LL^2(\GG))\bar\otimes \N$: $\gamma=(\id\otimes\alpha)(\cdot)_{[213]}$ is simply the amplified action, and $\beta=\Ad(\vv^{\GG *}_{[21]})\gamma$ is its cocycle twist. By the biduality theorem \cite[Theorem 2.6]{VaesUnitary}, there is an isomorphism $\Psi \colon \B(\LL^2(\GG))\bar\otimes \N\rightarrow \whG^{\op}\ltimes(\GG\ltimes\N)$ such that
\[
\alpha^{\ltimes \ltimes}\Psi(z)=(j_{\GG}R_{\GG}\otimes\Psi)\beta(z)\quad(z\in \B(\LL^2(\GG))\bar\otimes \N),
\]
where $j_{\GG}=J_{\vp_{\GG}}(\cdot)^* J_{\vp_{\GG}}$ is an $*$-anti-isomorphism $\LL^{\infty}(\GG)\rightarrow \LL^{\infty}(\GG)'$. Since $j_{\GG} R_{\GG}\colon \LL^{\infty}(\GG)\rightarrow \LL^{\infty}({\GG^{\op}}')$ is an isomorphism of quantum groups, weight $\wh{\wh\theta}\Psi$ on $\B(\LL^2(\GG))\bar\otimes \N$ is relatively $\delta_{\GG}^{-1}$-invariant under $\beta$. 
By Proposition \ref{prop1}, the standard implementation of $\beta$ with respect to $\wh{\wh{\theta}}\Psi$ satisfies the universal positive cone condition for $\omega$.

By Lemma \ref{lemma6} the same is true for $\mrm{U}^\beta$, the standard implementation of $\beta$ with respect to weight $\Tr\otimes \theta\in \oon{W}_{\nsf}(\B(\LL^2(\GG))\bar\otimes \N)$.\\

(Step 2) Let $\mrm{U}^\gamma$ be the standard implementation of $\gamma$ with respect to $\Tr\otimes \theta$. Specify the GNS Hilbert space for $\Tr\otimes \theta$ as $\msf{H}_{\Tr\otimes \theta}=\LL^2(\GG)\otimes \ov{\LL^2(\GG)}\otimes \msf{H}_\theta$, which comes from the identification $\HS(\LL^2(\GG))=\LL^2(\GG)\otimes \ov{\LL^2(\GG)}$. In this representation we have $\mrm{U}^\gamma = \mrm{U}^{\alpha}_{[14]}\in \LL^{\infty}(\GG)\bar\otimes \B(\LL^2(\GG)\otimes\ov{\LL^2(\GG)}\otimes \msf{H}_\theta)$ (Example \ref{ex4}) and \cite[Proposition 4.2]{VaesUnitary}
\begin{equation}\label{eq35}
\mrm{U}^\beta=(R_{\GG}\otimes j_{\LL^2(\GG)})(\vv^{\GG}_{[21]})_{[13]} \mrm{U}^\alpha_{[14]} \vv^{\GG}_{[21]}
\end{equation}
where $j_{\LL^2(\GG)}\colon \B(\LL^2(\GG))\rightarrow \B(\ov{\LL^2(\GG)})\colon j_{\LL^2(\GG)}(x)\ov{\xi}=\ov{x^* \xi}$. After applying $\Delta_{\GG}^{u,r}$ to the first leg of \eqref{eq35}, after small calculation we arrive at its universal version
\[
\mathds{U}^\beta=(R^u_{\GG}\otimes j_{\LL^2(\GG)})(\vV^{\GG}_{[21]})_{[13]} \mathds{U}^\alpha_{[14]} \vV^{\GG}_{[21]}
\]
and taking adjoint gives
\begin{equation}\label{eq36}
\mathds{U}^{\beta*}=
\vV^{\GG *}_{[21]}
 \mathds{U}^{\alpha *}_{[14]} 
(R^u_{\GG}\otimes j_{\LL^2(\GG)})(\vV^{\GG *}_{[21]})_{[13]}.
\end{equation}Applying $\pi$ to the first leg of \eqref{eq36} followed by vector functional, we obtain
\begin{equation}\label{eq37}\begin{split}
&\quad\;
\bigl\la 
\chi_1\otimes \eta_1\otimes \ov{\rho_1}\otimes \kappa_1\big|(\pi\otimes\id^{\otimes 3})(\mathds{U}^{\beta *})(
\chi_2\otimes \eta_2\otimes \ov{\rho_2}\otimes \kappa_2)\bigr\ra \\
&=\!
\bigl\la \chi_1 \!\otimes \kappa_1 \big|
\bigl(\pi\bigl(\!
(\omega_{\eta_1,\eta_2}\otimes \id)(\vV^{\GG *})\bigr)\!\otimes \!\I\bigr)
(\pi\otimes\id)(\mathds{U}^{\alpha *})\!
\bigl(\pi\bigl(
(\omega_{\ov{\rho_1},\ov{\rho_2}}j_{\LL^2(\GG)}\otimes R^u_{\GG})(\vV^{\GG *})\bigr)\!\otimes\! \I\bigr)
(\chi_2\otimes \kappa_2)\bigr\ra \\
&=
\bigl\la
\pi\bigl(
(\omega_{\eta_2,\eta_1}\otimes \id)(\vV^{\GG })\bigr)
 \chi_1\otimes \kappa_1 \big|
(\pi\otimes\id)(\mathds{U}^{\alpha *})\bigl(\pi\bigl(
(\omega_{\rho_2,\rho_1}R_{\widecheck{\GG}}\otimes \id )(\vV^{\GG *})\bigr)\chi_2\otimes \kappa_2\bigr)\bigr\ra \\
&=
\bigl\la
\pi\bigl(
(\omega_{\eta_2,\eta_1}\otimes \id)(\vV^{\GG })\bigr)
 \chi_1\otimes \kappa_1 \big|
(\pi\otimes\id)(\mathds{U}^{\alpha *})\bigl(\pi\bigl(
(\omega_{J_{\vp_{\GG}}\rho_1,J_{\vp_{\GG}}\rho_2}\otimes \id )(\vV^{\GG *})\bigr)\chi_2\otimes \kappa_2\bigr)\bigr\ra 
\end{split}\end{equation}

for all $\chi_i\in \msf{K}, \eta_i,\rho_i\in \LL^2(\GG),\kappa_i\in \msf{H}_\theta\,(i\in\{1,2\})$. Take $(\zeta_{i,j})_{i,j=1}^{n}, (\zeta'_{i,j})_{i,j=1}^{n}\in \mf{P}_{\Tr_n\otimes \theta}$. To prove the theorem, we need to show that
\begin{equation}\label{eq38}
\bigl\la 
(\zeta_{i,j}')_{i,j=1}^{n} \big| 
((\omega_{\xi_i,\xi_j}\pi\tau^u_{i/4}\otimes\id)(\mathds{U}^{\alpha *})\zeta_{i,j})_{i,j=1}^{n} \bigr\ra \ge 0.
\end{equation}
Assume first that 
\begin{equation}\label{eq40}
\xi_i=\pi \tau^{\GG, u}_{-i/4}\bigl((\omega_{\eta_i,\eta'_i}\otimes \id)(\vV^{\GG})\bigr)\tilde{\xi}_i
\end{equation}
where $\eta_i,\eta'_i\in \LL^2(\GG)$, $\tilde{\xi}_i\in \bigcup_{\eps>0}1_{\left[\eps,\eps^{-1}\right]}(A)\msf{K}$. Note that $(\omega_{\eta_i,\eta'_i}\otimes \id)(\vV^{\GG})\in \Dom(S_{\GG}^u)=\Dom(\tau^{\GG, u}_{-i/2})$ $\subseteq \Dom(\tau^{\GG,u}_{-i/4})$. Then $\xi_i\in \Dom(A^{\pm 1/4})$,
\[
A^{- 1/4} \xi_i=
A^{- 1/4} \pi \tau^{\GG,u}_{-i/4}\bigl(
(\omega_{\eta_i,\eta'_i}\otimes\id)(\vV^{\GG})\bigr)\tilde{\xi}_i=
\pi \bigl(
(\omega_{\eta_i,\eta'_i}\otimes\id)(\vV^{\GG})\bigr)A^{-1/4}\tilde{\xi}_i
\]
and
\[\begin{split}
&\quad\;
A^{ 1/4} \xi_j=
A^{ 1/4} \pi \tau^{\GG,u}_{-i/4}\bigl(
(\omega_{\eta_j,\eta'_j}\otimes\id)(\vV^{\GG})\bigr)\tilde{\xi}_j=
\pi \tau^{\GG,u}_{-i/2}\bigl(
(\omega_{\eta_j,\eta'_j}\otimes\id)(\vV^{\GG})\bigr)A^{ 1/4}\tilde{\xi}_j\\
&=
\pi R_{\GG}^u \bigl(
(\omega_{\eta_j,\eta'_j}\otimes\id)(\vV^{\GG *})\bigr)A^{ 1/4}\tilde{\xi}_j=
\pi \bigl(
(\omega_{\eta_j,\eta'_j}R_{\wc{\GG}}\otimes\id)(\vV^{\GG *})\bigr)A^{ 1/4}\tilde{\xi}_j\\
&=
\pi \bigl(
(\omega_{J_{\vp_{\GG}}\eta'_j,J_{\vp_{\GG}}\eta_j}\otimes\id)(\vV^{\GG *})\bigr)A^{ 1/4}\tilde{\xi}_j,
\end{split}\]

using $S^u_{\GG}=R^u_{\GG}\tau^{\GG,u}_{-i/2}$ and $(R_{\widecheck{\GG}} \otimes R^u_{\GG} )(\vV^{\GG})=\vV^{\GG}$. With these properties and \eqref{eq37}, we can write \eqref{eq38} as
\[\begin{split}
&\quad\;
\bigl\la 
(\zeta_{i,j}')_{i,j=1}^{n} \big| 
((\omega_{\xi_i,\xi_j}\pi\tau^u_{i/4}\otimes\id)(\mathds{U}^{\alpha *})\zeta_{i,j})_{i,j=1}^{n} \bigr\ra =
\sum_{i,j=1}^{n}
\bigl\la 
 \zeta_{i,j}' \big| 
(\omega_{A^{-1/4}\xi_i,A^{1/4}\xi_j}\pi \otimes\id)(\mathds{U}^{\alpha *})\zeta_{i,j} \bigr\ra \\
&=
\!\sum_{i,j=1}^{n}\!\!
\bigl\la 
\pi((\omega_{\eta_i,\eta'_i}\!\otimes\id)(\vV^{\GG}))A^{-1/4}\tilde{\xi}_i\!\otimes  \zeta_{i,j}' \big| \!
(\pi\!\otimes\id)(\mathds{U}^{\alpha *})\!
\bigl(
\pi((\omega_{J_{\vp_{\GG}} \eta_j',J_{\vp_{\GG}}\eta_j}\otimes\id)(\vV^{\GG *}))A^{1/4}\tilde{\xi}_j
\otimes \zeta_{i,j} \bigr)\bigr\ra \\
&=
\sum_{i,j=1}^{n}
\bigl\la 
A^{-1/4}\tilde{\xi}_i\otimes \eta'_i\otimes \ov{\eta'_j }\otimes \zeta'_{i,j} \big|
(\pi\otimes\id^{\otimes 3})(\mathds{U}^{\beta *})
(A^{1/4}\tilde{\xi}_j\otimes 
\eta_i\otimes \ov{ \eta_j}\otimes 
\zeta_{i,j})\bigr\ra\\
&=
\sum_{i,j=1}^{n}
\bigl\la 
 \eta'_i\otimes \ov{\eta'_j}\otimes \zeta'_{i,j}  \big|
(\omega_{A^{-1/4}\tilde{\xi_i},A^{1/4}\tilde{\xi}_j}\pi\otimes\id^{\otimes 3})(\mathds{U}^{\beta *})
(\eta_i\otimes \ov{\eta_j}\otimes \zeta_{i,j})\bigr\ra \\
&=
\bigl\la 
\bigl(\eta'_i\otimes \ov{\eta'_j}\otimes \zeta'_{i,j}  \bigr)_{i,j=1}^{n}\big|
\bigl(
(\omega_{ \tilde{\xi_i}, \tilde{\xi}_j}\pi\tau^{\GG, u}_{i/4}\otimes\id^{\otimes 3})(\mathds{U}^{\beta *})
(\eta_i\otimes \ov{\eta_j}\otimes \zeta_{i,j})\bigr)_{i,j=1}^{n}
\bigr\ra .
\end{split}\]
This expression is $\ge 0$ since
\begin{equation}\label{eq39}
\bigl(\eta'_i\otimes \ov{\eta'_j}\otimes \zeta'_{i,j}  \bigr)_{i,j=1}^{n}, \;
\bigl(
 \eta_i\otimes \ov{\eta_j}\otimes \zeta_{i,j}\bigr)_{i,j=1}^{n}
\in \mf{P}_{\Tr_n\otimes \Tr\otimes \theta}.
\end{equation}
Indeed, (up to taking sums and limits) we have $\zeta_{i,j}=x_i J_\theta \Lambda_\theta(x_j)$ for some $x_i\in \mf{N}_\theta$ (Lemma \ref{lemma3}) and then
\[\begin{split}
&\quad\;
\eta_i\otimes \ov{\eta_j}\otimes x_iJ_\theta \Lambda_\theta(x_j)=
\bigl( |\eta_i\ra\la \Omega|\otimes \I\otimes x_i\bigr) (\Omega\otimes \ov{\eta_j}\otimes J_\theta \Lambda_\theta(x_j) )\\
&=
\bigl( |\eta_i\ra\la \Omega|\otimes \I\otimes x_i\bigr) J_{\Tr\otimes \theta} (\eta_j\otimes \ov{\Omega}\otimes \Lambda_\theta(x_j) )\\
&=
\bigl( (|\eta_i\ra\la \Omega|\otimes \I)\otimes x_i\bigr) J_{\Tr\otimes \theta}
\Lambda_{\Tr\otimes \theta}\bigl(
(|\eta_j\ra\la \Omega|\otimes \I)\otimes x_j
\bigr)
\end{split}\]
for an arbitrary unit vector $\Omega\in \LL^2(\GG)$; thus \eqref{eq39} holds by Lemma \ref{lemma3}.\\

The above reasoning establishes inequality \eqref{eq38}, under the assumption that vectors $\xi_1,\dotsc,\xi_n\in \msf{K}$ have special form. Take now general $\xi_1,\dotsc,\xi_n\in \Dom(A^{1/4})\cap\Dom(A^{-1/4})$. Choose net $(\omega_\lambda)_{\lambda\in \Lambda}$ in $\LL^1(\wc\GG)$ such that $((\omega_\lambda\otimes\id)(\vV^{\GG}))_{\lambda\in \Lambda}$ is contractive and converges to $\I$ in the strict topology of $\M(\mrm{C}_0^u(\GG))$. Define $\wt{\omega}_\lambda=\tfrac{1}{\sqrt{\pi}}\int_{\RR} e^{-t^2} \omega_\lambda \tau^{\wc\GG}_{t}\md t$, then $\pi(a_\lambda)\xrightarrow[ \lambda\in\Lambda]{\sot}\I$, where $a_\lambda=(\wt{\omega}_\lambda\otimes\id)(\vV^{\GG})=\tfrac{1}{\sqrt{\pi}} \int_{\RR} e^{-t^2} \tau^{u,\GG}_t \bigl((\omega_\lambda\otimes\id)\vV^{\GG}\bigr)\md t$. Write $\tilde{\xi}_{i,m}=1_{\left[m^{-1},m\right]}(A)\xi_i \,(m\in\NN)$, and since $\LL^{\infty}(\wc\GG)\subseteq \B(\LL^2(\GG))$ is standard, we can write $\wt{\omega}_\lambda=\omega_{\eta_\lambda,\eta'_\lambda}$ for vectors $\eta_\lambda,\eta'_\lambda\in\LL^2(\GG)$. For each $1\le i\le n$, vector
\[
\xi_{i,\lambda,m}=\pi \tau^{\GG,u}_{-i/4}\bigl((\omega_{\eta_{\lambda},\eta'_\lambda}\otimes\id)(\vV^{\GG})\bigr)\tilde{\xi}_{i,m}
\]
 is of the form \eqref{eq40} and
 \[\begin{split}
&\quad\;
A^{\pm 1/4} \xi_{i,\lambda,m}=
A^{\pm 1/4} \pi \tau^{\GG,u}_{-i/4}\bigl(\tfrac{1}{\sqrt{\pi}}\int_{\RR} e^{-t^2} 
( \omega_{\lambda}\tau^{\wc\GG}_t \otimes\id)(\vV^{\GG})\md t\bigr)1_{\left[m^{-1},m\right]}(A)\xi_{i}\\
&=
\pi \tau^{\GG,u}_{\mp i/4 -i/4}\bigl(\tfrac{1}{\sqrt{\pi}}\int_{\RR} e^{-t^2} 
\tau^{\GG,u}_t (( \omega_{\lambda}  \otimes\id)(\vV^{\GG}))\md t\bigr)1_{\left[m^{-1},m\right]}(A)A^{\pm 1/4}\xi_{i}\\
&=
\tfrac{1}{\sqrt{\pi}}\int_{\RR} e^{-(t\pm i/4+i/4)^2} 
A^{it} \pi \bigl(
 ( \omega_{\lambda}  \otimes\id)(\vV^{\GG})  \bigr)A^{-it} 1_{\left[m^{-1},m\right]}(A)A^{\pm 1/4} \xi_{i} \md t
 \xrightarrow[(\lambda,m)\in \Lambda\times \NN]{}A^{\pm 1/4}\xi_i.
\end{split} \]
With this, we obtain
\[\begin{split}
&\quad\;
\bigl\la 
(\zeta_{i,j}')_{i,j=1}^{n} \big| 
((\omega_{\xi_i,\xi_j}\pi\tau^u_{i/4}\otimes\id)(\mathds{U}^{\alpha *})\zeta_{i,j})_{i,j=1}^{n} \bigr\ra \\
&=
\bigl\la 
(\zeta_{i,j}')_{i,j=1}^{n} \big| 
((\omega_{A^{-1/4} \xi_i,A^{1/4}\xi_j}\pi \otimes\id)(\mathds{U}^{\alpha *})\zeta_{i,j})_{i,j=1}^{n} \bigr\ra \\
&=
\lim_{(\lambda,m)\in \Lambda\times \NN}
\bigl\la 
(\zeta_{i,j}')_{i,j=1}^{n} \big| 
((\omega_{A^{-1/4} \xi_{i,\lambda,m},A^{1/4}\xi_{j,\lambda,m}}\pi \otimes\id)(\mathds{U}^{\alpha *})\zeta_{i,j})_{i,j=1}^{n} \bigr\ra \ge 0
\end{split}\]
which establishes \eqref{eq38} in general.
\end{proof}

As a corollary of Theorem \ref{thm5} we show that the operation of taking the standard implementation (with fixed weight) respects morphisms of quantum groups; in particular restriction of a standard implementation to a closed quantum subgroup is a standard implementation.

To set the stage, let $\GG,\HH$ be locally compact quantum groups with morphism $\GG\rightarrow \HH$; it is described by a non-degenerate $*$-homomorphism $\Pi\colon \mrm{C}_0^u(\HH)\rightarrow \M(\mrm{C}_0^u(\GG))$ which respects coproducts $(\Pi\otimes\Pi)\Delta^u_{\HH}=\Delta^u_{\GG}\Pi$ \cite{DKSS_ClosedSub, homomorphisms}. We can extend $\Pi$ to a normal, unital $*$-homomorphism $\Pi^{\vN}\colon \mrm{C}_0^u(\HH)^{**}\rightarrow \mrm{C}_0^u(\GG)^{**}$. Similarly, the reducing map $\pi_{\GG}\colon \mrm{C}_0^u(\GG)\rightarrow \mrm{C}_0(\GG)$ can be extended to a normal, unital $*$-homomorphism $\pi_{\GG}^{\vN} \colon \mrm{C}_0^u(\GG)^{**}\rightarrow \LL^{\infty}(\GG)$. Assume now that we have a von Neumann algebra $\N$ and action $\alpha_{\HH}\colon \HH\curvearrowright \N$. It admits the universal lift $\alpha^{u}_{\HH}\colon \N\rightarrow \mrm{C}_0^u(\HH)^{**}\bar\otimes \N$ (see Preliminaries), which can be used to define an action of $\GG$, $\alpha_{\GG}=(\pi_{\GG}^{\vN}\Pi^{\vN}\otimes \id)\alpha_{\HH}^{u}\colon \GG\curvearrowright \N$. One can also characterise $\alpha_{\GG}$ without using bidual \cst-algebras. See e.g.~\cite[Section 2.2]{DeCommerKrajczok} and \cite[Section 6.5]{DeCommerPhD} in the case of closed quantum subgroups, but the construction works for arbitrary morphisms. Compare also \cite[Section 2.3]{Kitamura} in \cst-setting. One should think of $\alpha_{\GG}$ as the composition $\GG\rightarrow \HH\curvearrowright \N$.

\begin{corollary}\label{cor1}
Let $\GG,\HH$ be locally compact quantum groups with morphism $\GG\rightarrow\HH$ described by $\Pi\colon \mrm{C}_0^u(\HH)\rightarrow \M(\mrm{C}_0^u(\GG))$ and $\N$ a von Neumann algebra with weight $\theta\in \oon{W}_{\nsf}(\N)$. Assume that we have actions $\alpha_{\HH}\colon \HH\curvearrowright \N, \alpha_{\GG}\colon \GG\curvearrowright \N$ as described in the previous paragraph. Denote by $\mrm{U}^{\alpha_{\HH}}, \mrm{U}^{\alpha_{\GG}}$ the standard implementations constructed with respect to $\theta$. Then
\[
\mathds{U}^{\alpha_{\GG}}=( \Pi\otimes\id)(\mathds{U}^{\alpha_{\HH}})\quad\textnormal{and}\qquad
\mrm{U}^{\alpha_{\GG}}=(\pi_{\GG} \Pi\otimes\id)(\mathds{U}^{\alpha_{\HH}}).
\]
\end{corollary}

\begin{proof}
It is enough to prove the second equation, the first follows by uniqueness of lift. To do this, we will establish that $(\pi_{\GG} \Pi\otimes\id)(\mathds{U}^{\alpha_{\HH}})$ satisfies conditions of Theorem \ref{thm1} for action $\alpha_{\GG}$. First, clearly it is a unitary representation of $\GG$ on $\msf{H}_\theta$. It implements $\alpha_{\GG}$: for $n\in\N$ we have
\[\begin{split}
(\pi_{\GG} \Pi\otimes\id)(\mathds{U}^{\alpha_{\HH}})^* &(\I\otimes n)
(\pi_{\GG} \Pi\otimes\id)(\mathds{U}^{\alpha_{\HH}})=
(\pi_{\GG}^{\vN}\Pi^{\vN}\otimes\id)\bigl(
\mathds{U}^{\alpha_{\HH} *}(\I\otimes n)
\mathds{U}^{\alpha_{\HH}}\bigr)\\
&=
(\pi_{\GG}^{\vN}\Pi^{\vN}\otimes\id)\alpha_{\HH}^{u}(n)=\alpha_{\GG}(n).
\end{split}\]
Next,
\[\begin{split}
&\quad\;
(J_{\vp_{\whG}}\otimes J_{\theta})
(\pi_{\GG} \Pi\otimes\id)(\mathds{U}^{\alpha_{\HH}})^*
(J_{\vp_{\whG}}\otimes J_{\theta})=
(R_{\GG}\pi_{\GG}\Pi\otimes j_\theta)(\mathds{U}^{\alpha_{\HH}})\\
&=
(\pi_{\GG}\Pi R_{\HH}^u \otimes j_\theta)(\mathds{U}^{\alpha_{\HH}})=
(\pi_{\GG}\Pi\otimes \id)(\mathds{U}^{\alpha_{\HH}})
\end{split}\]
where $j_{\theta}=J_{\theta}(\cdot)^* J_{\theta}$. In this calculation we have used  Proposition \ref{prop3} and the fact that strong quantum homomorphisms respect unitary antipodes, a consequence of \cite[Proposition 3.10]{homomorphisms}. Finally, take $n\in\NN, \xi_1,\dotsc,\xi_n\in \Dom(P_{\GG}^{-1/4})\cap\Dom(P_{\GG}^{1/4})$. Then
\[
\bigl(
(\omega_{\xi_i,\xi_j}\tau^{\GG}_{i/4}\otimes \id)\bigl( (\pi_{\GG}\Pi\otimes \id)\mathds{U}^{\alpha_{\HH} *}\bigr) \bigr)_{i,j=1}^{n}=
\bigl(
(\omega_{\xi_i,\xi_j} \pi_{\GG}\Pi \tau^{u, \HH}_{i/4} \otimes \id)(\mathds{U}^{\alpha_{\HH} *}) \bigr)_{i,j=1}^{n}
\]
preserves the positive cone by Theorem \ref{thm5}: as the representation in Definition \ref{def2} we take $\pi_{\GG}\Pi\colon \mrm{C}_0^u(\HH)\rightarrow \B(\LL^2(\GG))$, Hilbert space $\msf{K}$ is $\LL^2(\GG)$ and for the operator $A$ we take $P_{\GG}$. 
\end{proof}

\bibliographystyle{plain}

\end{document}